\documentclass[a4paper,12pt]{amsart}

\usepackage{amsmath,times,epsfig,amssymb,amsbsy,amscd,amsfonts,amstext,color,bm}

\usepackage[all]{xy}
\usepackage{bm}
\usepackage{comment}
\usepackage{enumerate}
\usepackage{circledsteps}
\usepackage{geometry}

\usepackage{tikz,epic,eso-pic}

\tikzset{v/.style={circle, draw, inner sep=2pt, minimum size=6pt, fill=white}}

\theoremstyle{plain}
\newtheorem{theorem}{Theorem}[section]
\newtheorem*{theorem*}{Theorem}
\newtheorem*{theoremA*}{Theorem A}
\newtheorem*{theoremB*}{Theorem B}
\newtheorem{corollary}[theorem]{Corollary}

\theoremstyle{definition}
\newtheorem{definition}[theorem]{Definition}
\newtheorem{remark}[theorem]{Remark}
\newtheorem{example}[theorem]{Example}
\newtheorem{proposition}[theorem]{Proposition}
\newtheorem{convention}[theorem]{Convention}

\newcommand{\Ker}{\operatorname{Ker}}

\newcommand{\Hom}{\operatorname{Hom}}

\newcommand{\corank}{\operatorname{corank}}
\newcommand{\coker}{\operatorname{Coker}}

\newcommand{\rank}{\operatorname{rank}}

\newcommand{\cl}{\operatorname{cl}}
\newcommand{\codim}{\operatorname{codim}}
\newcommand{\Char}{\operatorname{char}}

\newcommand{\bC}{\mathbb{C}}

\newcommand{\bK}{\mathbb{K}}

\newcommand{\bZ}{\mathbb{Z}}

\newcommand{\cA}{\mathcal{A}}

\newcommand{\cC}{\mathcal{C}}

\newcommand{\cL}{\mathcal{L}}

\newcommand{\cP}{\mathcal{P}}

\newcommand{\cR}{\mathcal{R}}

\title{The cohomology ring of the boundary manifold of a combinatorial line arrangement}
\author{Sakumi Sugawara}
\date{\today}
\address{Department of Mathematics, Faculty of Science, Hokkaido University, North 10, West 8, Kita-ku, Sapporo 060-0810, JAPAN. }
\email{math.skm.sugawara@gmail.com}

\begin{document}
\maketitle

\begin{abstract}
We compute the integral cohomology rings of the associated $3$-manifold called the \textit{boundary manifold} for a combinatorial line arrangement, and prove that it is isomorphic to the double of the Orlik-Solomon algebra. 
This generalizes the Cohen-Suciu doubling formula for complex realizable case to arbitrary combinatorial line arrangements, including non-realizable ones.
To handle the non-realizable case without geometric realization, we construct explicit homology cycles and compute intersection products, following the method of Doig--Horn for graph manifolds.  
As an application, we derive several results on the resonance variety of the boundary manifold.
\end{abstract}


\section{Introduction}
A finite set of lines in a projective plane is called a \textit{line arrangement}.
Line arrangements are studied from various viewpoints, including combinatorics, algebraic geometry, and topology, as fundamental objects in the theory of hyperplane arrangements \cite{orl-ter}.
From a topological viewpoint, a central question is whether various topological invariants are combinatorially determined, that is, determined by the data of intersection lattices. For instance, the cohomology ring of the complement is combinatorially determined and is known as the Orlik-Solomon algebra, whereas the fundamental group is not.
It remains an open problem whether many invariants of the complement 
are combinatorially determined.

The boundary of the regular neighborhood of the union of lines (called the \textit{boundary manifold}) has also been studied \cite{coh-suc-bound}.
It is known that the boundary manifold is a closed $3$-manifold which is obtained from the incidence graph of the line arrangement by plumbing along the graph \cite{JY}. 
Hence, the homeomorphism type of the boundary manifold is combinatorially determined. 
Boundary manifolds have also been used to investigate the topology of the complement \cite{fgm, hiro}, for example.

There is also a combinatorial generalization of line arrangements, called a \textit{combinatorial line arrangement}, which abstracts the projective geometric property that any two lines intersect at a single point (also called line combinatorics; see \cite{accm}, for example).
As in the case of ordinary line arrangements, a combinatorial line arrangement determines an incidence graph, which is useful for studying the combinatorial property.

Recently, Ruberman--Starkston studied the topological realizability of combinatorial line arrangements \cite{rub-sta}.
They not only studied realizability in various categories (such as topological and symplectic) but also constructed, from a combinatorial line arrangement, a closed $3$-manifold that coincides with the boundary manifold in realizable cases, and investigated their contact-geometric properties.
Their construction of the $3$-manifold is based on plumbing along the incidence graph.
In this paper, we study the topology of such $3$-manifolds, which we call the \textit{boundary manifold of a combinatorial line arrangement}.
Concerning the cohomology ring of the boundary manifold of a complex projective line arrangement (=realizable cases), Cohen--Suciu proved the following result.

\begin{theorem}\label{thm:coh-suc} (\cite{coh-suc-proj}, Theorem 4.2) \,
Let $\cA$ be a line arrangement in $\bC P^2$, X be its complement, and $M$ be the boundary manifold. Then, the cohomology ring $H^{*}(M; \bZ)$ is isomorphic to the double of the cohomology ring of the complement as graded algebras.
\end{theorem}

The definition of the \textit{double} for a graded algebra will be given in Section 2.3.
As mentioned above, since the topological type of $M$ and the cohomology ring of $X$ are described combinatorially,
the statement of the theorem is combinatorial.
Their proof relies on the long exact sequence in cohomology groups and on Hodge-theoretic property, and is therefore not ``combinatorial'' (however, the above theorem is a special case of a more general result valid for general hypersurfaces).

The main result of this paper is a purely combinatorial proof of the above doubling formula with its generalization to arbitrary combinatorial line arrangements. 
\begin{theorem}
Let $\cC$ be a combinatorial line arrangement and $M$ be its boundary manifold. Then, the cohomology ring $H^{*}(M; \bZ)$ is isomorphic to the double of the Orlik-Solomon algebra of $\cC$.
\end{theorem}

Since the proof of Theorem \ref{thm:coh-suc} relies on a complex realization, we have to take a different strategy valid for arbitrary combinatorial line arrangements.
Our approach for the proof is to compute the intersection ring, which is dual to the cohomology ring of $M$.
Specifically, we construct a basis of the homology groups and observe how their intersections behave.
The method for computing the intersection ring of a plumbed $3$-manifold is developed by Doig--Horn \cite{doi-hor}, and we follow their approach. 
Finally, we deduce that the cohomology ring is isomorphic to the double of the Orlik-Solomon algebra, which is defined for any combinatorial line arrangement.
Note that, for realizable cases, this gives an alternative proof of Theorem \ref{thm:coh-suc}.
Doig--Horn computed the intersection ring over rational coefficients for general plumbed 3-manifolds. However, we computed the product structure over \textit{intergral} coefficients for the boundary manifolds. 

As an application of our main result, we give some results on the resonance varieties of the cohomology ring of the boundary manifold.
The resonance variety for the doubled algebra was also studied in \cite{coh-suc-proj}. The techniques developed in that paper are all applicable to our setting, and the results concerning the resonance variety can also be generalized to combinatorial line arrangements.

This paper is organized as follows. Sections 2 and 3 provide the necessary preliminaries for stating the main theorem. In particular, we review plumbed manifolds in Section 3.
Section 4 is devoted to defining the boundary manifold and constructing a basis of the homology groups. This construction is one of the most important parts of the paper.
In Section 5, we compute the intersection product using the basis constructed in Section 4.
Finally, we compute the resonance variety of the boundary manifold in Section 6.

\vspace{3mm}
\textbf{Acknowledgements.}
The author would like to thank Yongqiang Liu, Takuya Saito, Alexandru Suciu and Masahiko Yoshinaga for valuable comments and suggestions on this paper.

\section{Preliminaries}
\subsection{Combinatorial line arrangements}

\begin{definition}
A \textit{combinatorial line arrangement} is a pair $\cC = (\cL, \cP)$ where $\cL$ is a finite set and $\cP \subset 2^{\cL}$ such that 
\begin{enumerate}[(i)]
\item For each $P \in \cP$, $|P| \geq 2$,
\item For distinct $\ell, \ell' \in \cL$, there exists a unique $P \in \cP$ such that $\ell, \ell' \in P$.
\end{enumerate}
An element in $\cL$ (resp. \! $\cP$) is called a \textit{line} (resp. \! \textit{intersection point}).
\end{definition}

For $\ell_i \in \cL$ and $P_I \in \cP$, let denote $\cP_{i} = \{P \in \cP \mid P \ni \ell_i\}$ the set of intersection points containing $\ell_i$ and $\cL_I = \{\ell \in \cL \mid \ell \in P_I\}$ the set of lines that contained in $P_I$.
For lines $\ell_j, \ell_k$, there exists a unique intersection point. Let $I(j,k)$ be the index set such that $\ell_j, \ell_k \in P_{I(j,k)}$.

\begin{example}
Let $\bK$ be a field and $\cA = \{\ell_0, \ell_1, \ldots, \ell_n\}$ be a line arrangement in the projective plane $\bK P^2$. Let $\cP(\cA)$ be the set of intersection points of lines in $\cA$. 
Since each intersection point is labelled by lines passing through it, we can consider $\cP (\cA)$ as a subset of $2^\cA$. 
Then, the pair $(\cA, \cP (\cA))$ is a combinatorial line arrangement.
\end{example}

\begin{definition}
Let $\bK$ be a field.
A combinatorial line arrangement $\cC=(\cL, \cP)$ is \textit{realizable over $\bK$} if there exists a line arrangement $\cA$ in $\bK P^2$ and bijections $\psi_1: \cL \rightarrow \cA$ and $\psi_2: \cP \rightarrow \cP (\cA)$ such that $\ell \in P$ if and only if $\psi_1 (\ell) \ni \psi_2 (P)$ for each $\ell \in \cL$ and $P \in \cP$.
This arrangement $\cA$ is called the \textit{realization} of $\cC$.
\end{definition}

We show an example of a combinatorial line arrangement that is not realizable over any field.

\begin{figure}[htbp]
\centering
\begin{tikzpicture}
\draw (0,1) --++(4,0.5);
\fill (0,1)++(0.5,0.0625) circle (0.06);
\fill (0,1)++(2,0.25) circle (0.06);
\fill (0,1)++(3.5,0.4375) circle (0.06);
\draw (0,-1) --++(4,-0.5);
\fill (0,-1)++(0.5,-0.0625) circle (0.06);
\fill (0,-1)++(2,-0.25) circle (0.06);
\fill (0,-1)++(3.5,-0.4375) circle (0.06);
\draw[red] (0.5, 1.0625) -- (2,-1.25);
\draw[red] (0.5, -1.0625) -- (2,1.25);
\draw[blue] (0.5, 1.0625) -- (3.5,-1.4375);
\draw[blue] (0.5, -1.0625) -- (3.5,1.4375);
\draw[green] (2, 1.25) -- (3.5,-1.4375);
\draw[green] (2, -1.25) -- (3.5,1.4375);

\draw[densely dashed] (0,0) --++(2.4,0) to[out=90,in=180] ++(0.3,0.3) to[out=0,in=90] ++(0.3,-0.3)--++(0.5,0); 

\filldraw (1.17,0) circle (0.06) node[below] {\footnotesize $P$};
\filldraw (1.8,0) circle (0.06) node[below] {\footnotesize $Q$};
\filldraw (2.7,0) circle (0.06) node[below] {\footnotesize $R$};

\end{tikzpicture}
\caption{A combinatorial line arrangement which is not realizable over any field by Pappus' theorem.}
\label{fig:Pappus}
\end{figure}
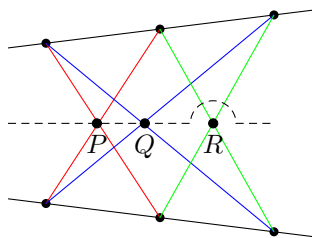

\begin{example}
There is a combinatorial line arrangement that cannot be realized over any field.
For the configuration in Figure \ref{fig:Pappus}, a line passing through $P$ and $Q$ must pass through $R$, by Pappus' theorem.
Therefore, the configuration containing a ``line'' passing through $P$ and $Q$ but not $R$ gives a non-realizable combinatorial line arrangement.
\end{example}

\begin{figure}[htbp]
\centering
\begin{tikzpicture}
\coordinate (P) at (0,-2);
\coordinate (Q) at (5,0);
\coordinate (R) at (5,-3.5);

\draw (P)++(0,-0.3) --++(0,2.3) node[above]{$\ell_1$};
\draw (P)++(-0.2,-0.2) --++(1.7,1.7) node[right]{$\ell_2$};
\draw (P)++(-0.3,0) --++(2.3,0) node[right]{$\ell_3$};
\draw (P)++(-0.3,2)--++(2.3,-2.3)node[below]{$\ell_{0}$} ;
\draw (P)++(2,-0.17) --++(-2.3,1.3) node[left]{$\ell_{4}$};

\draw (P) ++ (1,-0.5) node[below]{$\cC$};

\draw (P) ++ (1, -2) node{$\textbf{nbc}(\cC) =\{(1,2),(1,3),(1,4),(2,4)\}$};
 
\draw (Q) node[below]{$\ell_{0}$};
\draw (Q)++(1,0) node[below]{$\ell_{1}$};
\draw (Q)++(2,0) node[below]{$\ell_{2}$};
\draw (Q)++(3,0) node[below]{$\ell_{3}$};
\draw (Q)++(4,0) node[below]{$\ell_{4}$};

\draw (Q)++(0,1.5) node[above]{$P_{01}$};
\draw (Q)++(1,1.5) node[above]{$P_{02}$};
\draw (Q)++(2,1.5) node[above]{$P_{034}$};
\draw (Q)++(3,1.5) node[above]{$P_{123}$};
\draw (Q)++(4,1.5) node[above]{$P_{14}$};
\draw (Q)++(5,1.5) node[above]{$P_{24}$};

\draw (Q) --++ (0,1.5);
\draw (Q) --++ (1,1.5);
\draw (Q) --++ (2,1.5);
\draw (Q)++(1,0) --++ (-1,1.5);
\draw (Q)++(1,0) --++ (2,1.5);
\draw (Q)++(1,0) --++ (3,1.5);
\draw (Q)++(2,0) --++ (-1,1.5);
\draw (Q)++(2,0) --++ (1,1.5);
\draw (Q)++(2,0) --++ (3,1.5);
\draw (Q)++(3,0) --++ (-1,1.5);
\draw (Q)++(3,0) --++ (0,1.5);
\draw (Q)++(4,0) --++ (-2,1.5);
\draw (Q)++(4,0) --++ (0,1.5);
\draw (Q)++(4,0) --++ (1,1.5);

\draw (Q) ++(2,-0.8) node{$\Gamma(\cC)$};

\draw (R) node[below]{$\ell_{0}$};
\draw (R)++(1,0) node[below]{$\ell_{1}$};
\draw (R)++(2,0) node[below]{$\ell_{2}$};
\draw (R)++(3,0) node[below]{$\ell_{3}$};
\draw (R)++(4,0) node[below]{$\ell_{4}$};

\draw (R)++(0,1.5) node[above]{$P_{01}$};
\draw (R)++(1,1.5) node[above]{$P_{02}$};
\draw (R)++(2,1.5) node[above]{$P_{034}$};
\draw (R)++(3,1.5) node[above]{$P_{123}$};
\draw (R)++(4,1.5) node[above]{$P_{14}$};
\draw (R)++(5,1.5) node[above]{$P_{24}$};

\draw (R) --++ (0,1.5);
\draw (R) --++ (1,1.5);
\draw (R) --++ (2,1.5);
\draw (R)++(1,0) --++ (-1,1.5);
\draw (R)++(1,0) --++ (2,1.5);
\draw (R)++(1,0) --++ (3,1.5);
\draw (R)++(2,0) --++ (-1,1.5);
\draw[dotted] (R)++(2,0) --++ (1,1.5);
\draw (R)++(2,0) --++ (3,1.5);
\draw (R)++(3,0) --++ (-1,1.5);
\draw[dotted] (R)++(3,0) --++ (0,1.5);
\draw (R)++(4,0) --++ (-2,1.5);
\draw[dotted] (R)++(4,0) --++ (0,1.5);
\draw[dotted] (R)++(4,0) --++ (1,1.5);

\draw (R) ++(2,-0.8) node{A spanning tree $T$};

\end{tikzpicture}
\caption{A combinatorial line arrangement and its incidence graph}
\label{fig:arr}
\end{figure}
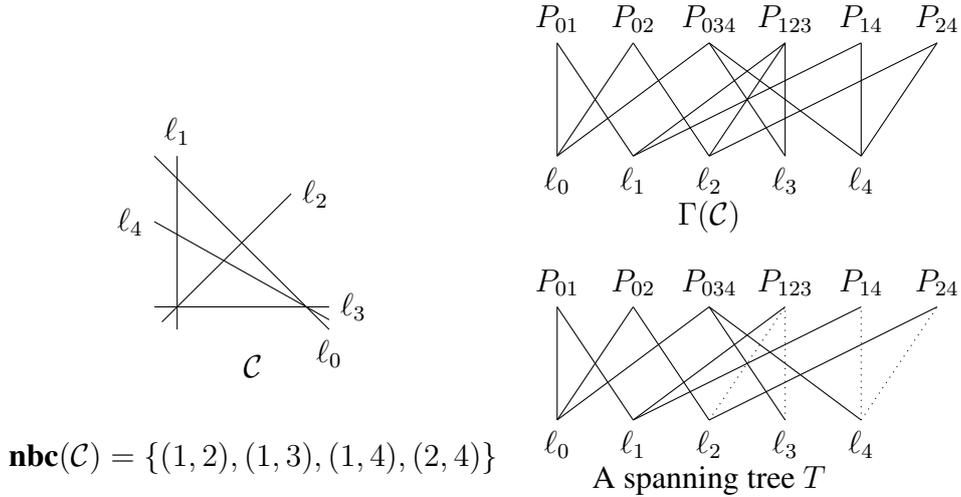

\begin{definition}
Let $\cC = (\cL, \cP)$ be a combinatorial line arrangement.
Define a bipartite graph $\Gamma = \Gamma(\cC)$ by the vertex set $V(\Gamma) = \cL \cup \cP$ and the edge set $E(\Gamma) = \{(\ell, P) \mid \ell \in P\}$.
The graph $\Gamma$ is called the \textit{incidence graph} of a combinatorial line arrangement (see Figure \ref{fig:arr}).
\end{definition}

\begin{remark}
If $\cC$ is realized as a line arrangement in a projective space, the graph $\Gamma(\cC)$ has the same information as the intersection poset of the line arrangement, which is a well-known object. 
In general, $\Gamma (\cC)$ has the same information as the lattice of flats of the simple matroid of rank $3$ obtained from the combinatorial line arrangement.
Lines and points of a combinatorial line arrangement correspond to rank $1$ and $2$ flats, respectively. 
The Hasse diagram of the lattice of flats is considered as a graph, and $\Gamma(\cC)$ is the graph obtained by removing the minimum and maximum elements from it.
\end{remark}

\begin{definition}
Let $\cC = (\cL, \cP)$ be a combinatorial line arrangement, and the lines are ordered as $\cL = \{\ell_0, \ell_1, \ldots, \ell_n\}$. 
We fix the first line $\ell_0$ as a distinguished one.
A pair $S = (j,k)$ of integers ($1 \leq j < k \leq n$) is called a \textit{nbc}, if there exists $P \in \cP$ such that $S \subset P$ and $\min S = \min P$.
Denote $\textbf{nbc}(\cC)$ the set of all nbc in $\cC$.
\end{definition}

Note that the index $0$ of the distinguished line does not appear in $S$.
For example, if $\cC = (\cL, \cP)$ is a pencil, that is, every line is contained in the same one point, then $\textbf{nbc}(\cC) = \emptyset$. It is because the unique intersection point $P$ satisfies $\min P = 0$ and for any $S = (j,k)$ with positive integers, $\min S \neq 0$.

\begin{remark}
The word \textit{nbc} originates from the term \textit{non broken circuit} in the theory of matroids and hyperplane arrangements, see \cite{bjo} or Sections 3.1, 3.2 in \cite{orl-ter}, for example.
If $\cC$ is realizable, $\textbf{nbc}(\cC)$ defined here corresponds to the nbc basis of rank $2$ of the Orlik-Solomon algebra for an affine line arrangement obtained by deconing the line $\ell_0$. 
\end{remark}

From now on, we will always assume that lines in a combinatorial line arrangement are ordered by $\cL = \{\ell_0, \ell_1, \ldots, \ell_n\}$ and $\ell_0$ is fixed as the distinguished line.

\begin{proposition}\label{prop:nbc}
There is a one-to-one correspondence between independent cycles of the incidence graph $\Gamma(\cC)$ and elements in $\textbf{nbc}(\cC)$.
\end{proposition}

\proof
At first, an independent cycle of a connected graph corresponds to an edge of the complement of a spanning tree.
Define a subset $T$ of the edge set of $\Gamma (\cC)$ by
\[
T = \{(P_{I}, \ell_{i}) \mid \ell_0 \in P_{I} \} \cup \{(P_{I}, \ell_{i})\mid i = \min I \}.
\]
Then, $T$ is a spanning tree and the complement $C = E \setminus T$ is written as 
\[
C = \{(P_I, \ell_{j}) \mid \ell_0 \notin P_{I} , j \neq \min I\}.
\]
To each edge in $(P_{I}, \ell_j)$ in $C$, we assign the nbc $(\min I, j)$.
Next, consider the converse. For a nbc $S=(i,j)$, the intersection point $P_{I(i,j)}$ of $\ell_i$ and $\ell_{j}$ satisfies $i= \min I(i,j)$. 
So, we assign the edge $(P_{I(i,j)}, \ell_{j})$ to this nbc.
\endproof

\subsection{Orlik-Solomon algebra}
Let $\cC = (\cL, \cP)$ be a combinatorial line arrangement with $\cL = \{\ell_0 ,\ell_1, \ldots, \ell_n\}$.
Define a graded $\bZ$-algebra $A = A(\cC) = \oplus_{k=0}^{2} A^{k}$ as follows: $A^0 =\bZ$, $A^1$ is a free $\bZ$-module with a basis $\{e_{1}, \ldots ,  e_{n} \}$, symbols corresponding to lines $\ell_{1}, \ldots, \ell_{n}$, and $A^{2}$ is a free $\bZ$-module with a basis $\{f_{j,k} \mid (j,k) \in \textbf{nbc} (\cC)\}$. 
The multiplication $\mu: A^1 \wedge A^1 \rightarrow A^2$ is defined by 
\begin{eqnarray*}
\mu (e_{i}, e_{j}) = 
\left \{
\begin{array}{ll}
f_{i,j} & \mbox{if $(i,j) \in \textbf{nbc}(\cC)$}, \\
f_{k,j}-f_{k,i} & \mbox{if there exists $k$ such that $(k,i),(k,j)\in \textbf{nbc}(\cC)$}, \\
0 & \mbox{otherwise}.
\end{array}
\right .
\end{eqnarray*}
for $1 \leq i < j \leq n$. This graded algebra $A$ is called the \textit{Orlik-Solomon algebra} of $\cC$.

\begin{remark}
If $\cC$ is realizable over some field, then $A(\cC)$ is isomorphic to the Orlik-Solomon algebra of an affine line arrangement obtained by deconing along the line $\ell_{0}$.
The isomorphism class of the Orlik-Solomon algebra does not depend on the choice of the distinguished line (see \cite{dim-hyp}, Theorem 3.4).
\end{remark}

\begin{remark}
Suppose that $\cC$ is realized as a line arrangement $\cA$ in $\bC P^2$. Then, the Orlik-Solomon algebra $A$ is isomorphic to the cohomology ring of the complement $X = \bC P^2 \setminus \bigcup_{H \in \cA}H$ \cite{orl-sol}. 
In particular, $b_2 (X) = |\textbf{nbc} (\cC)|$.
\end{remark}

\subsection{Doubling graded algebras}
We will explain the doubling of a graded algebra in a general setting here. For the details, see Section 3.1 of \cite{coh-suc-proj}.
Let $R$ be a commutative ring and $A= \oplus_{k=0}^{m} A^k$ be a connected, finitely generated graded $R$-algebra.
The dual $\overline{A} = \Hom_{R}(A, R)$ is also a graded $R$-algebra with the grading $\overline{A}^k = \Hom_R (A^k, R)$.
For $a \in A$ and $f \in \overline{A}$, by defining $(a \cdot f) (b) = f(ba)$ and $(f \cdot a) (b) = f(ab)$, $\overline{A}$ is a $A$-bimodule. 
Note that if $a \in A^{k}$ and $f \in \overline{A}^j$, then $af, fa \in \overline{A}^{j-k}$.

We define the \textit{double} $\widehat{A} = \textsf{D}(A)$ of a graded algebra $A$ by the base $A$-module $\widehat{A}= A \oplus \overline{A}$ with the multplication $(a,f)\cdot (b,g) = (ab, ag + bf)$. 
The double is also graded $R$-algebra, with the grading $\widehat{A}^{k}=A^{k} \oplus \overline{A}^{2m-1-k}$.

In this paper, we deal with two-dimensional graded algebras.
Let $m=2$ and assume that each $A^k$ is a free module. First, we have 
\begin{eqnarray*}
\widehat{A} &=& \widehat{A}^{0} \oplus \widehat{A}^{1} \oplus \widehat{A}^{2} \oplus \widehat{A}^{3} \\
&=& A^0 \oplus (A^1 \oplus \overline{A}^2) \oplus (A^2 \oplus \overline{A}^1) \oplus \overline{A}^0
\end{eqnarray*}
Let $\mu: A^1 \otimes A^1 \rightarrow A^2$ be the multipication on $A$ and write $\mu(a_i,a_j) = \sum_{k=1}^{r_2} \mu_{i,j,k}b_{k}$ with specific basis 
$\{a_i\}_{i=1}^{r_{1}}$ and $\{b_{k}\}_{k=1}^{r_{2}}$, where $r_{p}$ is the rank of $A^{p}$ for $p=1,2$.
Then, the multipilication of $\widehat{\mu}: \widehat{A}^1 \otimes \widehat{A}^1 \rightarrow \widehat{A}^2$ is the same as $\mu$ on $A^1 \oplus A^1$, vanishes on $\overline{A}^2 \oplus \overline{A}^2$, and 
\[
\widehat{\mu} (a_{j}, \overline{b}_{k}) = \sum_{i=1}^{r_1} \mu_{i,j,k} \overline{a}_{i}
\]
on $A^1 \oplus \overline{A}^2$, where $\{ \overline{a}_{i} \}_{i=1}^{r_{1}}$ and $\{ \overline{b}_k \}_{k=1}^{r_{2}}$ are the dual basis.

\begin{example}
Let $A=A(\cC)$ be the Orlik-Solomon algebra of a combinatorial line arrangement. 
Note that $\rank (A^{1} \oplus \overline{A}^{2}) = n+|\textbf{nbc}(\cC)|= n+b_{1} (\Gamma (\cC))$, from Proposition \ref{prop:nbc}. The multiplication on $A^{1} \oplus \overline{A}^2$ of the double is computed as follows (see also Section 7 in \cite{coh-suc-bound}):
\begin{eqnarray*}
\widehat{\mu}(e_{i}, \overline{f}_{j,k})= 
\left \{
\begin{array}{ll}
-\overline{e}_{i} + \sum_{m \in I(j,i)} \overline{e}_{m} & (i=k), \\
-\overline{e}_{k} & (i \in I(j,k), i \neq k), \\
0 & (\mbox{otherwise}).
\end{array}
\right .
\end{eqnarray*}
Note that if $\cC$ is a pencil, then $\widehat{\mu}$ vanishes on $\widehat{A}^{1} \otimes \widehat{A}^{1}$.
\end{example}

\section{Plumbed manifolds}
In this section, we recall some basic matters on plumbed manifolds, which will be needed in this paper. For details, see \cite{neu}, for example.

\subsection{$S^1$-bundles over $S^2$}
Firstly, we recall the construction of oriented $S^1$-bundles.
Let $\pi: E \rightarrow B$ be an oriented $S^1$-bundle.
As is well known, the isomorphism class of oriented $S^1$-bundles is classified by elements in $H^2 (B; \bZ)$, which is called the Euler class.
In particular, if $B$ is an orientable surface, then $H^2 (B; \bZ) \cong \bZ$.
The isomorphism class of oriented $S^1$-bundles over $B$ is classified by an integer, which is called the \textit{Euler number} of the oriented $S^1$-bundle.

For the case where $B=S^2$, let us describe a specific construction of $S^1$-bundle $\pi: E \rightarrow S^2$ with a given Euler number.
Let $D^{triv} \subset S^2$ be a closed disk and $S'= \cl(S^2 \setminus D^{triv})$.
Since both $D^{triv}$ and $S'$ are homeomorphic to a disk, $S^1$-bundles restricted over $D^{triv}$ and $S'$ are trivial, so we describe them as $\pi^{-1} (D^{triv}) \cong D^{triv} \times S^1$ and $\pi^{-1} (S') \cong S' \times S^1$.
Let $\varphi: \partial D^{triv} \times S^1 \rightarrow \partial S' \times S^1$ be the gluing map.
By giving orientations in $D^{triv}, S'$ and the fiber $S^1$, the gluing map is represented by a quadratic square matrix of integer coefficients with the determinant $-1$. In particular, since the gluing map preserves the fiber, the matrix is described by a matrix 
$ \bigl(
\begin{smallmatrix}
   -1 & e \\
   0 & 1
\end{smallmatrix}
\bigl)$.
The integer $e$ that appeared here is nothing but the Euler number of this $S^1$-bundle.

\subsection{Plumbed manifolds}
Next, we will explain plumbed manifolds, which are obtained by gluing $S^1$-bundles over surfaces.
Let $\Gamma = (V, E)$ be a connected, simple finite graph and suppose that there is a vertex weight $e: V \rightarrow \bZ; i \mapsto e(i)$ over integers.
In this paper, such a graph is called a \textit{plumbing graph}.
From a plumbing graph, we can define an oriented closed $3$-manifold as follows.
First, we assign an $S^1$-bundle $\pi_{i}: E_{i} \rightarrow S_{i}$ over $S_{i} \cong S^2$ with the Euler number $e(i)$ to each vertex $i$ of $\Gamma$.
As discussed above, decompose the base space $S_{i}$ into two disks $D_{i}^{triv}$ and $S'_{i}$ and trivialize the bundles over them.
Let $V_{i}$ be the set of vertices adjacent to $i$.
For each vertex $j \in V_{i}$, take a disjoint disk $D_{i}^{j}$ in $S'_{i}$. Let denote $S''_{i} = S'_{i} \setminus \bigcup_{j \in V_{i}} D_{i}^{j}$ the complement.
For a vertex $j \in V_{i}$, glue a part of the boundaries of $\partial D_{i}^{j} \times S^1 (\subset \pi_{i}^{-1} (S''_{i}) )$ and $\partial D_{j}^{i} \times S^1 (\subset \pi^{-1}_{j} (S''_{j}))$ by a homeomorphism represent by the matrix
$ \bigl(
\begin{smallmatrix}
   0 & 1 \\
   1 & 0
\end{smallmatrix}
\bigl)$.
By repeating this operation for each edge of $\Gamma$, we obtain an oriented closed manifold $M(\Gamma)$. The resulting manifold is called a \textit{plumbed manifold}.
The disks $D_{i}^{j}$ and the tori $\partial D_{i}^{j} \times \partial D_{j}^{i}$ used in the gluing are called plumbing disks and plumbing tori, respectively.

\begin{remark}
In general, the term plumbed manifolds refers to a larger class that may involve surfaces of non-zero genus or gluings with negative signs \cite{neu}.
However, only a restricted subclass of plumbed manifolds defined here appears in this paper. So we will restrict our attention to this case.
\end{remark}

\subsection{The first homology group of a plumbed manifold}
For a plumbing graph, define a integer matrix $A(\Gamma) = (a_{i,j})_{i,j \in V}$ as follows:
\begin{eqnarray*}
a_{i,j} = 
\left \{
\begin{array}{ll}
e(i) & (i=j), \\
1 & (\{i,j\} \in E), \\
0 & (\mbox{otherwise}).
\end{array}
\right .
\end{eqnarray*}
Using this matrix, we can compute the first homology group of a plumbed manifold (see Lemma 3.4 and the proof of Theorem 3.2 in \cite{acm}, or Theorem 3.6 in \cite{doi-hor}, for example).
\begin{proposition}
\[
H_{1} (M(\Gamma); \bZ) \cong \bZ^{b_{1} (\Gamma)} \oplus \coker (A(\Gamma)).
\]
In particular, $b_{1} (M(\Gamma)) = b_{1} (\Gamma) + \corank (A(\Gamma))$ and the torsion part is computed by the Smith normal form of $A(\Gamma)$.
\end{proposition}

\begin{remark}
Here, we mention the explicit generators of $H_{1} (M(\Gamma); \bZ)$ (see Section 3 in \cite{doi-hor}). 
The summand $\bZ^{b_{1} (\Gamma)}$ is generated by $1$-cycles in $M(\Gamma)$ obtained from cycles of the graph $\Gamma$.
On the other hand, the matrix $A(\Gamma)$ is considered as an endomorphism of a free $\bZ$-module $\bZ^{|V|}$ generated by $1$-cycles represented by the fiber $S^1$ of the bundle corresponding to each vertex.
Thus, $\coker(A(\Gamma))$ is generated by the fiber $S^1$ corresponding to some vertices.
\end{remark}

\section{The boundary manifold}
In this section, we define the boundary manifold of a combinatorial line arrangement and study its topology.
\subsection{The definition}
Again, let $\cC = (\cL, \cP)$ be a combinatorial line arrangement ordered as $\cL = \{\ell_{0}, \ell_{1}, \ldots, \ell_{n}\}$ with the distinguished line $\ell_{0}$.
For the incidence graph $\Gamma = \Gamma(\cC)$, define the vertex weight $e: V(\Gamma) \rightarrow \bZ$ by $e(\ell_{i}) = 1 - |\cP_{i}|$ for each line and $e(P_{I}) = -1$ for each intersection point. Recall that $\cP_{i} = \{P \in \cP \mid P \ni \ell_{i} \}$. 
With this weight, we consider $\Gamma(\cC)$ as a plumbing graph.
The plumbed manifold $M=M(\cC)$ obtained from the plumbing graph $\Gamma(\cC)$ is called the \textit{boundary manifold} of a combinatorial line arrangement $\cC$.

\begin{remark}
If $\cC$ is realizable over $\bC$ with a realization $\cA = \{\ell_0, \ell_1 , \ldots, \ell_{n}\}$, then the manifold $M$ is diffeomorphic to the boundary of a regular neighborhood of $\bigcup_{i=0}^{n} \ell_i$. This is why we call the manifold $M$ the \textit{boundary manifold} in this paper.
See \cite{coh-suc-bound, hiro, JY} for details of this object.
\end{remark}

\begin{remark}
Even for the case where $\cC$ is not realizable over $\bC$, the manifold $M$ is studied by Ruberman--Starkston \cite{rub-sta} with a symplectic/contact geometric motivation.
They call the manifold $M$ \textit{the canonical contact manifold for the line arrangement} (with a certain contact structure) in that paper. They study the tightness of the contact structure and the relationship between symplectic fillability and symplectic realizability of $M$.
\end{remark}

\begin{remark}
The plumbing graph $\Gamma(\cC)$ is not a ``reduced'' one. 
In fact, we can blow down the vertices corresponding to intersection points consisting of two lines, and obtain a plumbing graph with fewer vertices.
Many previous studies have used this reduced graph as the plumbing graph. However, we use the plumbing graph defined here to address all intersection points in the same way.
\end{remark}

\subsection{Notations and conventions}
We will introduce and summarise the symbols that will be needed later on. 
At first, the notations regarding the vertex $\ell_{i}$:
\begin{itemize}
\item $\pi_{i}: E_{i} \rightarrow S_{i}$: the $S^1$-bundle over $S_{i} \cong S^2$ with the Euler number $1 -|\cP_{i}|$, corresponding to the vertex $\ell_{i}$ ($i=0,1, \ldots, n$).
\item $D^{triv}_{i} (\subset S_{i})$: the trivializing disk.
\item $D^{I}_{i} (\subset S_{i}) $: the plumbing disk corresponding to the edge $(P_{I}, \ell_{i})$.
\item $S_{i}'' = S_{i} \setminus (D^{triv}_{i} \cup \bigcup_{P_{I} \in \cP_{i}} D^{I}_{i} )$.
\item $\alpha(i;I_1, I_2)$: a simple proper arc in $S_{i}''$ connecting $\partial D_{i}^{I_1}$ and $\partial D_{i}^{I_2}$ ($P_{I_1}, P_{I_2} \in \cP_i$).
\item $A(i; I_1, I_2) = \pi^{-1}_{i} (\alpha(i;I_1, I_2))$: the fiber annulus.
\item $\partial_{I_{k}} A(i;I_1, I_2) = A(i;I_1, I_2) \cap \pi_{i}^{-1} (\partial D^{I_{k}}_{i}) $: one of the connected components of the boundary of the annulus ($k=1,2$).
\item $\alpha(i;triv, I)$: a simple proper arc in $S_{i}''$ connecting $\partial D^{triv}_{i}$ and $\partial D^{I}_{i}$ ($P_{I} \in \cP_{i}$).
\item $A(i;triv, I) = \pi^{-1}_{i}(\alpha(i;triv, I))$: the fiber annulus.
\item $\partial_{I} A(i;triv, I) = A(i;triv, I) \cap \pi^{-1}_{i}(\partial D^{I}_{i})$: a connected component of the boundary of the annulus.
\end{itemize}
Similarly, we define the notations regarding the vertex $P_I$:
\begin{itemize}
\item $\pi_{I}: E_{I} \rightarrow S_{I}$: the $S^1$-bundle over $S_{I} \cong S^2$ with the Euler number $-1$, corresponding to the vertex $P_{I}$ ($P_I \in \cP$).
\item $D^{triv}_{I} (\subset S_{I})$: the trivializing disk.
\item $D^{i}_{I} (\subset S_{I}) $: the plumbing disk corresponding to the edge $(P_{I}, \ell_{i})$.
\item $S_{I}'' = S_{I} \setminus (D^{triv}_{I} \cup \bigcup_{\ell_{i} \in \cL_{I}} D^{i}_{I})$.
\item $\alpha(I;i_1, i_2)$: a simple proper arc in $S_{I}''$ connecting $\partial D_{I}^{i_1}$ and $\partial D_{I}^{i_2}$ ($\ell_{i_1}, \ell_{i_2} \in \cL_I$).
\item $A(I; i_1, i_2) = \pi^{-1}_{I} (\alpha(I;i_1, i_2))$: the fiber annulus.
\item $\partial_{k} A(I;i_1, i_2) = A(I;i_1, i_2) \cap \partial \pi^{-1}_{I} (D^{i_{k}}_{I})$: one of the connected components of the boundary of the annulus ($k=1,2$).
\item $\alpha(I;triv, i)$: a simple proper arc in $S_{I}''$ connecting $\partial D^{triv}_{I}$ and $\partial D^{i}_{I}$ ($\ell_{i} \in \cP_{I}$).
\item $A(I;triv, i) = \pi^{-1}_{I}(\alpha(I;triv, i))$: the fiber annulus.
\item $\partial_{i} A(I;triv, i) = A(I; triv, i) \cap \pi^{-1}_{I}( \partial D^{i}_{I} )$: a connected component of the boundary of the annulus.
\end{itemize}

We assume the following for the arcs:
\begin{itemize}
\item[(A1)] For $j=1,\ldots, n$, the arc $\alpha(j; I(j,0), I(j,k_1))$ does not intersect with $\alpha (j; I(j,0), I(j,k_2))$ for each $k_1, k_2$ with $I(j, k_1) \neq I(j,k_2)$.
\item[(A2)] The arc $\alpha (0;triv,I(i,0))$ does not intersect with $\alpha (0;I(j,0), I(k,0))$ for each $i,j,k$.
\item[(A3)] For $P_{I}\in \cP$, the arc $\alpha (I;triv, i)$ does not intersect $\alpha(I; j,k)$ for each $i,j,k$ with $\ell_i, \ell_j, \ell_k \in P_{I}$.
\item [(A4)] For $I \in \cP_{0}$, the arc $\alpha (I;i,0)$ does not intersect with $\alpha(I; j,0 ) $ for distinct $\ell_i, \ell_j \in P_{I}$.
\end{itemize}

\begin{convention}\label{conv:orientation}
Here, we remark on the orientation of submanifolds, following Guillemin--Pollack \cite{gul-pol}.
Let $M$ be an oriented manifold and $X$, $Y$ be submanifolds of $M$.
We assume that $X$ and $Y$ intersect transversely and $\dim X + \dim Y > \dim M$.
We consider the orientation of the transversal intersection $X \cap Y$.
Let $x \in X\cap Y$. We orient the normal verctor spaces $N_{x} (X|M)$ and $N_{x}(Y|M)$ so that 
\begin{eqnarray*}
N_{x}(X|M) \oplus T_{x} X = T_{x}M \\
N_{x}(Y|M) \oplus T_{x} Y = T_{x}M 
\end{eqnarray*}
as oriented vector spaces. Then, we give the orientation to $T_{x}(X \cap Y)$ so that 
\[
N_{x}(X|M) \oplus N_{x}(Y|M) \oplus T_{x} (X \cap Y) = T_{x}M 
\]
as an oriented vector space. Note that $X \cap Y = (-1)^{(\codim X)(\codim Y)} (Y \cap X)$.

Let us give an orientation to the boundary manifold $M=M(\cC)$.
Since $M$ is obtained from $S^1$-bundles over surfaces by gluing each other,
it suffices to orient a piece of $S^1$-bundles.
Fix a line $\ell_i$ and consider the piece $\pi_i^{-1}(S''_{i}) \cong S_{i}'' \times S^1$.
We fix orientations to $S''_{i}$ and $S^1$ and orient $S_{i}'' \times S^1$ in this order.
The boundary $\partial (S''_{i} \times S^1 )$ is oriented so that the outward normal vector is the first summand in the tangent space.
The manifold $M$ is oriented by extending the orientation defined here.

\end{convention}

\subsection{The first homology group of $M$}
\begin{proposition}\label{prop:1sthomology}
$H_{1} (M; \bZ)$ is a free $\bZ$-module of rank $b_{1} (\Gamma(\cC)) + n$.
\end{proposition}

\proof

It is sufficient to compute the cokernel of the matrix $A(\Gamma)$.
By observing each row, we have only the following two types of equations as relations in the homology group:
\begin{eqnarray*}
(1-|\cP_{i}|) t_{i}  + 
\sum_{I, P_{I} \ni \ell_{i}} t_{I} = 0 & \, (0 \leq i \leq n), \\ 
-t_{I} + \sum_{j, \ell_{j} \in P_{I}} t_{j} = 0 & \, (P_{I} \in \cP),
\end{eqnarray*}
where each $t_{i}, t_{I}$ is represented by the fiber of the $S^1$-bundle corresponding to $\ell_i, P_{I}$. 
From the second equation, $t_{I}$ is redundant for each $I$. By substituting the second equations into the first equations for each $I$, we have the following only one relation:
\[
\sum_{i=0}^{n} t_{i} = 0,
\]
for each $0\leq i \leq n$.
Thus, $t_{0} = - (t_{1} + \cdots +t_{n})$. 
Therefore, we have that $\coker (A(\Gamma))$ is a free $\bZ$-module generated by $\{t_{1} , \ldots, t_{n}\}$.
\endproof

Next, we will observe the explicit cycles forming a basis of $H_{1} (M; \bZ)$.
First, by the above computation, the fibers $S^1$ of the bundles corresponding to vertices $\ell_{1},  \ldots, \ell_{n}$ represent a part of the basis.
Fix a point $x_{i} \in S_{i}''$ ($i=1, \ldots, n$) and
let $t_{i} = \{x_{i}\} \times S^1$ be the representing loop.
Such a cycle obtained from the fiber of each vertex is called \textit{$t$-type cycle}. 
Each cycle $t_{i}$ is positively oriented with respect to the fixed orientation.

$1$-cycles corresponding to cycles of the graph $\Gamma(\cC)$ are described as follows. 
First, recall that independent cycles of $\Gamma(\cC)$ correspond one-to-one to elements in $\textbf{nbc}(\cC)$ (Proposition \ref{prop:nbc}).
The cycle in $\Gamma(\cC)$ corresponding to $(j,k) \in \textbf{nbc}(\cC)$ is defined by the cycle concatenating $\ell_{0}, P_{I(k,0)}, \ell_{k}, P_{I(j,k)}, \ell_{j}, P_{I(j,0)}$ and $\ell_{0}$ in this order.
The corresponding cycle in $M$ is constructed by connecting simple arcs that connect the plumbing disks in the base spaces, following the cycle of this graph.
Explicitly, the oriented cycle concatenating $\alpha(0;I(j,0),I(k,0))$, $\alpha(I(k,0); 0,k)$, $\alpha(k; I(k,0), I(j,k))$, $\alpha(I(j,k);k,j)$, $\alpha(j; I(j,k),I(j,0))$ and $\alpha(I(j,0);j,0)$ in this order corresponds.
Let denote $\gamma_{j,k}$ this resulting cycle in $H_{1}(M; \bZ)$.
Such a cycle obtained from the cycle of $\Gamma(\cC)$ is called \textit{$\gamma$-type cycle}. 
Note that the tangent vector of $\gamma_{j,k}$ points outward from $S_{k}''$ on $\partial D_{k}^{I(j,k)}$.

\begin{remark}\label{rmk:realizable1cycle}
Suppose that $\cC$ is realizable over $\bC$ with a realization $\cA = \{\ell_{0}, \ell_{1}, \ldots, \ell_{n}\}$. Then a $t$-type cycle corresponds to the meridian of each line $\ell_{i}$ ($i=1, \ldots, n$) and a $\gamma$-type cycle corresponds to the cycle obtained from the triangle with vertices the intersections $\ell_0 \cap \ell_{j}$, $\ell_{0} \cap \ell_{k}$ and $\ell_{j} \cap \ell_{k}$ by pushing out to the boundary of the neighborhood $((j,k) \in \textbf{nbc}(\cC) )$.
\end{remark}

\subsection{The second homology group of $M$}
Next, we construct explicit $2$-cycles (embedded surfaces) representing a basis of $H_{2} (M; \bZ)$.
Since $M$ is an oriented closed $3$-manifold, each element in $H_{2} (M; \bZ)$ is described by the interesection dual of an element in $H_{1} (M; \bZ)$, by Poincar\'e duality.
That is, when $\{\alpha_{1}, \ldots, \alpha_{b_{1}(M)}\}$ represents a basis of $H_{1} (M; \bZ)$, the set $\{ID(\alpha_{1}) ,\ldots, ID(\alpha_{b_{1}(M)})\}$ of $2$-cycles satisfying $ID(\alpha_{i}) \cap \alpha_{j} = \delta_{i,j}$ represents the basis of $H_{2} (M; \bZ)$, where $\cap$ denotes the intersection product.
Now, we have two kinds of representatives of the basis of $H_{1} (M; \bZ)$, $t$-type and $\gamma$-type.
Thus, we will construct the intersection duals for each type.
We will follow the construction developed by Doig--Horn \cite{doi-hor} here.

\subsubsection{The dual of $\gamma$-type cycles}
Let $\gamma_{j,k} \in H_{1} (M; \bZ)$ be the $\gamma$-type cycle corresponding to $(j,k) \in \textbf{nbc}(\cC)$.
Let $\tau_{j,k} = \partial D^{I(j,k)}_{k} \times \partial D^{k}_{I(j,k)}$ be the plumbing torus corresponding to the edge $(P_{I(j,k)}, \ell_{k})$.
By the definition of the cycle $\gamma_{j,k}$, $\tau_{j,k}$ intersects with $\gamma_{j,k}$ at once, and does not intersect with other cycles of $\gamma$-type. 
Give an orientation to $\tau_{j,k}$ such that $\tau_{j,k} \cap \gamma_{j,k}$ is a positive intersection.
Since $t$-type cycles are represented by the fiber of each $S^1$-bundle piece, the intersection number of $t$-type cycles and $\gamma$-type cycles is $0$.
Therefore, the plumbing torus $\tau_{j,k}$ is the intersection dual of the $\gamma$-type cycle $\gamma_{j,k}$.
Note that the normal vector of $\tau_{j,k}$ points outward from $S_{k}''$ at each point (see Convention \ref{conv:orientation}).

\subsubsection{The dual of $t$-type cycles}
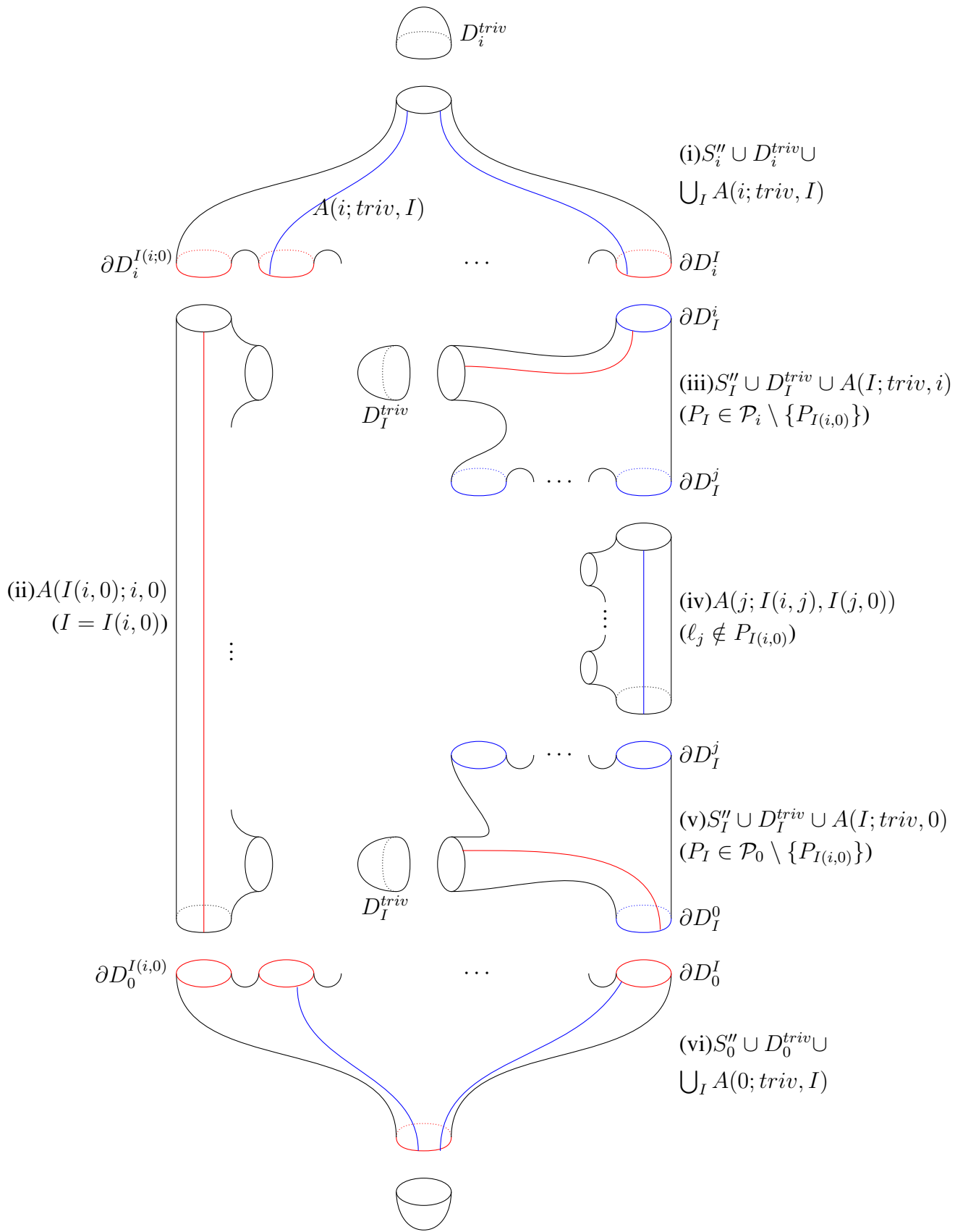
\begin{figure}[htbp]
\centering
\begin{tikzpicture}
\coordinate (P) at (0,0);
\coordinate (Q) at (4,-4);
\coordinate (R) at (4,-8);
\coordinate (S) at (4,-12);
\coordinate (T) at (-4,-4);
\coordinate (U) at (0,-19);

\draw (P) circle [x radius=0.5, y radius=0.25];
\draw (P) ++ (-0.5,0) to[out=-90,in=90] ++(-4,-3);
\draw[red] (P) ++(-4.5,-3) to[out=-90,in=180]++(0.5,-0.25) to[out=0,in=-90]++(0.5,0.25);
\draw (P)++ (-3.5,-3)to[out=90,in=180]++(0.25,0.25)to[out=0,in=90]++(0.25,-0.25);
\draw [red] (P)++(-3,-3)to[out=-90,in=180]++(0.5,-0.25)to[out=0,in=-90]++(0.5,0.25);
\draw (P)++(-2,-3)to[out=90,in=180]++(0.25,0.25)to[out=0,in=90]++(0.25,-0.25);
\draw (P) ++ (0.5,0) to[out=-90,in=90] ++(4,-3);
\draw[red] (P) ++ (4.5,-3)to[out=-90,in=0]++(-0.5,-0.25)to[out=180,in=-90]++(-0.5,0.25);
\draw (P) ++ (3.5,-3)to[out=90,in=0]++(-0.25,0.25) to[out=180, in=90]++(-0.25,-0.25);
\draw[densely dotted,red] (P) ++(-3.5,-3)to[out=90,in=0]++(-0.5,0.25)to[out=180,in=90]++(-0.5,-0.25);
\draw[densely dotted,red] (P) ++(-2,-3)to[out=90,in=0]++(-0.5,0.25)to[out=180,in=90]++(-0.5,-0.25);
\draw[densely dotted,red] (P) ++(4.5,-3)to[out=90,in=0]++(-0.5,0.25)to[out=180,in=90]++(-0.5,-0.25);
\draw (P) ++ (1,-3) node{$\cdots$};
\draw[densely dotted] (P) ++(0.5,1)to[out=90,in=0]++(-0.5,0.25)to[out=180,in=90]++(-0.5,-0.25);
\draw(P) ++(0.5,1)to[out=-90,in=0]++(-0.5,-0.25)to[out=180,in=-90]++(-0.5,0.25);
\draw (P) ++ (0.5,1)to[out=90,in=0]++(-0.5,0.7)to[out=180,in=90]++(-0.5,-0.7);
\draw[blue] (P) ++ (-0.3,-0.2)to[out=-90,in=90] ++ (-2.5,-3);
\draw[blue] (P) ++ (0.3,-0.2)to[out=-90,in=90] ++ (3.4,-3);
\draw (P)++(0.5,1.2)node[right]{$D^{triv}_{i}$};
\draw (P) ++ (4.5,-1)node[right]{(i)$S''_i \cup D_{i}^{triv} \cup$};
\draw (P)++(4.5,-1.7)node[right]{ $\bigcup_{I} A(i;triv,I)$ };
\draw (P)++(-4.4,-3)node[left]{$\partial D_{i}^{I(i;0)}$};
\draw (P)++(-1,-2)node{$A(i;triv,I)$};
\draw (P)++(4.5,-3)node[right]{$\partial D_{i}^{I}$};

\draw[blue] (Q) circle[x radius = 0.5, y radius=0.25];
\draw (Q) ++ (0.5,0) --++(0,-3);
\draw [blue] (Q) ++(0.5,-3) to[out=-90,in=0]++(-0.5,-0.25)to[out=180,in=-90]++(-0.5,0.25);
\draw (Q) ++(-0.5,-3)to[out=90,in=0]++(-0.25,0.25)to[out=180,in=90]++(-0.25,-0.25);
\draw (Q) ++ (-0.5,0) to[out=-90, in=0] ++(-3,-0.5);
\draw (Q) ++ (-0.5,0) ++ (-3,-0.5)++(0,-0.5) circle[x radius=0.25, y radius =0.5];
\draw (Q) ++ (-0.5,0) ++ (-3,-0.5)++(0,-1) to[out=0,in=90]++(1,-0.5)to[out=-90,in=90]++(-1,-1);
\draw[blue](Q)++(-3.5,-3) to[out=-90,in=180]++(0.5,-0.25)to[out=0,in=-90]++(0.5,0.25);
\draw (Q)++(-2.5,-3)to[out=90,in=180]++(0.25,0.25)to[out=0,in=90]++(0.25,-0.25);
\draw[densely dotted, blue] (Q) ++(0.5,-3)to[out=90,in=0]++(-0.5,0.25)to[out=180,in=90]++(-0.5,-0.25);
\draw[densely dotted, blue] (Q) ++(-2.5,-3)to[out=90,in=0]++(-0.5,0.25)to[out=180,in=90]++(-0.5,-0.25);
\draw (Q) ++ (-1.5,-3) node{$\cdots$};
\draw (Q) ++(-4.5,-1.5) to[out=0,in=-90]++(0.25,0.5)to[out=90,in=0]++(-0.25,0.5);
\draw[densely dotted] (Q) ++(-4.5,-1.5) to[out=180,in=-90]++(-0.25,0.5)to[out=90,in=180]++(0.25,0.5);
\draw (Q) ++ (-4.5,-1.5) to[out=180,in=-90]++(-0.7,0.5)to[out=90,in=180]++(0.7,0.5);
\draw[red] (Q)++(-0.2,-0.23)to[out=-90,in=0] ++(-3.05,-0.65);
\draw (Q) ++(0.5,-1.2) node[right]{(iii)$S''_{I}\cup D^{triv}_{I} \cup A(I;triv,i)$};
\draw (Q)++ (0.5,-1.8) node[right]{($P_I \in \cP_{i} \setminus \{P_{I(i,0)}\}$)};
\draw (Q) ++ (0.5,0) node[right]{$\partial D_{I}^{i}$};
\draw (Q) ++ (0.5,-3) node[right]{$\partial D_{I}^{j}$};
\draw (Q) ++(-4.7,-1.8)node{$D_{I}^{triv}$};

\draw (R) circle[x radius = 0.5, y radius=0.25];
\draw (R) ++(0.5,0)--++(0,-3) to[out=-90,in=0]++(-0.5,-0.25)to[out=180,in=-90]++(-0.5,0.25)to[out=90,in=0] ++(-0.5,0.3);
\draw (R) ++ (-0.5,0) to[out=-90,in=0] ++(-0.5,-0.3);
\draw (R) ++ (-0.5,0) ++(-0.5,-0.3) ++ (0,-0.3)circle[x radius = 0.15, y radius=0.3];
\draw (R) ++ (-0.5,-3) ++(-0.5,0.3) ++ (0,0.3)circle[x radius = 0.15, y radius=0.3];
\draw (R) ++ (-0.5,0) ++(-0.5,-0.3) ++ (0,-0.3) ++ (0,-0.3) to[out=0,in=90]++(0.3,-0.3);
\draw (R) ++ (-0.5,-3) ++(-0.5,0.3) ++ (0,0.3) ++ (0,0.3) to[out=0,in=-90]++(0.3,0.3);
\draw[densely dotted] (R) ++(0.5,-3)to[out=90,in=0]++(-0.5,0.25)to[out=180,in=90]++(-0.5,-0.25);
\draw (R) ++ (-0.7,-1.4) node{$\vdots$};
\draw[blue] (R) ++ (0,-0.25) --++(0,-3);
\draw (R) ++(0.5,-1.2) node[right]{(iv)$A(j;I(i,j),I(j,0))$};
\draw (R) ++(0.5,-1.8) node[right]{($\ell_{j} \notin P_{I(i,0)}$)};

\draw[blue] (S) circle[x radius = 0.5, y radius=0.25];
\draw[blue] (S)++(-3,0) circle[x radius = 0.5, y radius=0.25];
\draw (S)++(-0.5,0)to[out=-90,in=0]++(-0.25,-0.25)to[out=180,in=-90]++(-0.25,0.25);
\draw (S)++(-2,0)to[out=-90,in=0]++(-0.25,-0.25)to[out=180,in=-90]++(-0.25,0.25);
\draw (S) ++ (-3.5,-2)  circle[x radius = 0.25, y radius=0.5];
\draw (S) ++ (-3.5,0) to[out=-90,in=0]++(0.5,-1.5) --++(-0.5,0);
\draw (S) ++ (-3.5,-2.5) to[out=0,in=90] ++ (3,-0.5);
\draw (S) ++ (0.5,0) --++(0,-3);
\draw[blue] (S) ++ (0.5,-3) to[out=-90,in=0]++(-0.5,-0.25)to[out=180,in=-90]++(-0.5,0.25);
\draw[densely dotted,blue] (S) ++(0.5,-3)to[out=90,in=0]++(-0.5,0.25)to[out=180,in=90]++(-0.5,-0.25);
\draw (S) ++ (-1.5,0) node{$\cdots$};
\draw (S) ++(-4.5,-2.5) to[out=0,in=-90]++(0.25,0.5)to[out=90,in=0]++(-0.25,0.5);
\draw[densely dotted] (S) ++(-4.5,-2.5) to[out=180,in=-90]++(-0.25,0.5)to[out=90,in=180]++(0.25,0.5);
\draw (S) ++ (-4.5,-2.5) to[out=180,in=-90]++(-0.7,0.5)to[out=90,in=180]++(0.7,0.5);
\draw[red] (S) ++(0.3,-3.2)to[out=90,in=0] ++(-3.6,1.45);
\draw (S) ++(0.5,-1.2) node[right]{(v)$S''_{I}\cup D^{triv}_{I} \cup A(I;triv,0)$};
\draw (S)++ (0.5,-1.8) node[right]{($P_I \in \cP_{0} \setminus \{P_{I(i,0)}\}$)};
\draw (S) ++ (0.5,0) node[right]{$\partial D_{I}^{j}$};
\draw (S) ++ (0.5,-3) node[right]{$\partial D_{I}^{0}$};
\draw (S) ++(-4.7,-2.8)node{$D_{I}^{triv}$};

\draw (T)circle[x radius=0.5,y radius =0.25];
\draw (T)++(-0.5,0) --++(0,-11)to[out=-90, in=180]++(0.5,-0.25)to[out=0, in=-90]++(0.5,0.25);
\draw (T) ++(0.5,0) to[out=-90,in=180] ++(0.5,-0.5);
\draw (T) ++ (1,-1) circle[x radius=0.25,y radius =0.5];
\draw (T) ++ (1,-1.5) to[out=180,in=90] ++ (-0.5,-0.5);
\draw (T) ++(0.5,-11) to[out=90,in=180] ++ (0.5,0.5);
\draw (T) ++ (1,-10) circle[x radius=0.25,y radius =0.5];
\draw (T) ++(1,-9.5) to[out=180,in=-90] ++ (-0.5,0.5);
\draw[densely dotted] (T) ++(0.5,-11)to[out=90,in=0]++(-0.5,0.25)to[out=180,in=90]++(-0.5,-0.25);
\draw (T) ++ (0.5,-6) node{$\vdots$};
\draw[red] (T) ++ (0,-0.25) --++(0,-11);
\draw (T)++(-0.5,-5) node[left]{(ii)$A(I(i,0); i,0)$};
\draw (T)++(-0.5,-5.6) node[left]{($I=I(i,0)$)};

\draw[red] (U) ++(0.5,0)to[out=-90,in=0]++(-0.5,-0.25)to[out=180,in=-90]++(-0.5,0.25);
\draw[densely dotted,red] (U) ++(0.5,0)to[out=90,in=0]++(-0.5,0.25)to[out=180,in=90]++(-0.5,-0.25);
\draw (U) ++ (-0.5,0) to[out=90,in=-90] ++(-4,3);
\draw (U) ++ (0.5,0) to[out=90,in=-90] ++(4,3);
\draw[red] (U) ++ (-4,3)circle[x radius=0.5,y radius=0.25];
\draw[red] (U) ++ (-2.5,3)circle[x radius=0.5,y radius=0.25];
\draw[red] (U) ++ (4,3)circle[x radius=0.5,y radius=0.25];
\draw (U) ++ (-3.5,3) to[out=-90,in=180] ++(0.25,-0.25) to[out=0,in=-90]++(0.25,0.25);
\draw (U) ++ (-2,3) to[out=-90,in=180] ++(0.25,-0.25) to[out=0,in=-90]++(0.25,0.25);
\draw (U) ++ (3,3) to[out=-90,in=180] ++(0.25,-0.25) to[out=0,in=-90]++(0.25,0.25);
\draw (U) ++ (1,3) node{$\cdots$};
\draw (U)++(0,-1) circle[x radius=0.5, y radius=0.25];
\draw (U) ++ (0.5,-1) to[out=-90, in=0]++(-0.5,-0.7) to[out=180, in=-90]++(-0.5,0.7);
\draw[blue] (U) ++ (3.6,2.85) to[out=-120,in=90]++(-3.3,-3.1);
\draw[blue] (U) ++ (-2.3,2.75) to[out=-90,in=90]++(2.2,-3);
\draw (U) ++ (4.5,1.7)node[right]{(vi)$S''_0 \cup D_{0}^{triv} \cup$};
\draw (U)++(4.5,1)node[right]{ $\bigcup_{I} A(0;triv,I)$ };
\draw (U) ++(4.5,3)node[right]{$\partial D_{0}^{I}$};
\draw (U) ++(-4.5,3)node[left]{$\partial D_{0}^{I(i,0)}$};

\end{tikzpicture}
\caption{The intersection dual of $t_{i}$}
\label{fig:dualoft}
\end{figure}

Let $t_{i} \in H_{1} (M; \bZ)$ be the $t$-type cycle corresponding to the vertex $\ell_{i}$ ($i=1, \ldots, n$).
The construction of the intersection dual of a $t$-type cycle is more complicated than that of $\gamma$-type. 
The intersection dual is constructed by gluing several pieces of the $S^1$-bundles corresponding to each vertex. 
The construction is divided into six steps. See Figure \ref{fig:dualoft} for the construction.

In the following construction, manifolds with corners appear, but we assume that they are appropriately smoothed.

\underline{(i) The vertex $\ell_{i}$}

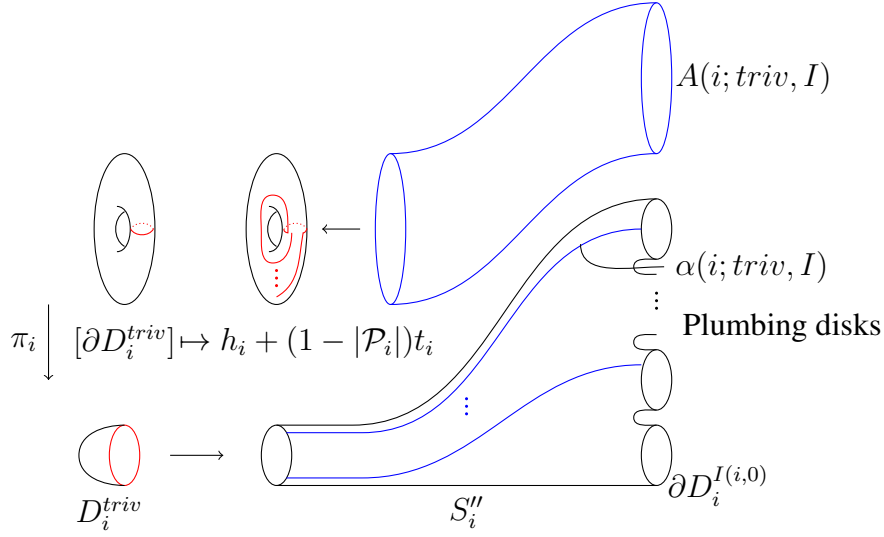
\begin{figure}[hbtp]
\centering
\begin{tikzpicture}
\coordinate (P) at (0,0);

\draw (P) circle[x radius =0.2, y radius=0.4];
\draw (P)++(0,-0.4) --++ (5,0);
\draw (P) ++ (5,0) circle[x radius =0.2, y radius=0.4];
\draw (P) ++ (5,1) circle[x radius =0.2, y radius=0.4];
\draw (P) ++ (5,3) circle[x radius =0.2, y radius=0.4];
\draw (P) ++ (0,0.4)--++(1,0) to[out=0, in=180]++ (4,3);
\draw (P) ++ (5,0.4) to[out=180,in=-90] ++(-0.3,0.1)to[out=90,in=180] ++(0.3,0.1);
\draw (P) ++ (5,1.4) to[out=180,in=-90] ++(-0.3,0.1)to[out=90,in=180] ++(0.3,0.1);
\draw (P) ++ (5,2.4) to[out=180,in=-90] ++(-0.3,0.1)to[out=90,in=180] ++(0.3,0.1);
\fill  (P)++(5,2.05)circle(0.02);
\fill  (P)++(5,2.15)circle(0.02);
\fill  (P)++(5,1.95)circle(0.02);

\draw (P) ++ (5.2,1.7)node[right]{Plumbing disks};
\draw (P) ++ (5,-0.4)node[right]{$\partial D_{i}^{I(i,0)}$};

\draw[blue] (P) ++ (0.13,0.3)--++(1,0) to[out=0,in=180] ++(3.67,2.7);
\draw[blue] (P) ++ (0.13,-0.3)--++(1,0) to[out=0,in=180] ++ (3.67,1.5);
\fill[blue]  (P)++(2.5,0.75)circle(0.02);
\fill[blue]  (P)++(2.5,0.55)circle(0.02);
\fill[blue]  (P)++(2.5,0.65)circle(0.02);

\draw (P) ++ (4,2.8) to[out=-90,in=180]++(1.1,-0.3)node[right]{$\alpha(i;triv,I)$};

\draw[->] (P) ++ (1.1,3) --++(-0.5,0);

\draw[blue] (P) ++ (1.5,3) circle [x radius=0.2, y radius=1 ];
\draw[blue] (P) ++ (1.5,3) ++(0,1) to[out=0,in=180] ++(3.5,2);
\draw[blue] (P) ++ (1.5,3) ++(0,-1) to[out=0,in=180] ++(3.5,2);
\draw [blue] (P) ++ (1.5,3) ++ (3.5,2) circle [x radius=0.2, y radius=1 ];
\draw (P) ++ (1.5,3) ++ (3.6,2) node [right] {$A(i;triv, I)$};

\draw (P) ++ (2.5,-0.4) node[below]{$S''_{i}$};

\draw[red] (P) ++(-0.5,0) ++ (-1.5,0) circle[x radius =0.2, y radius=0.4];
\draw (P) ++(-0.5,0)++ (-1.5,0.4) to[out=180,in=90] ++(-0.6,-0.4)to[out=-90,in=180] ++(0.6,-0.4);

\draw (P) ++(-0.5,0) ++(-1.7,-0.4)node[below]{$D_{i}^{triv}$};

\draw (P) ++ (0,3) circle[x radius =0.4, y radius=1];
\draw (P) ++ (0,3) ++ (-0.1,0.3) to[out=0, in=0] ++(0,-0.6);
\draw (P) ++ (0,3) ++ (0,0.23) to[out=210,in=150] ++ (0,-0.46);

\draw (P) ++ (0.4,1.5) node{$\mapsto h_{i} + (1-|\cP_{i}|)t_{i}$};

\draw[red] (P) ++ (0.1,3) to[out=-120, in =180] ++ (0.05,-0.05)to[out=90,in=0] ++ (-0.2,0.5)to[out=180,in=90]++(-0.15,-0.5)to[out=-90,in=180] ++ (0.15,-0.4)to[out=0, in =-90] ++(0.25,0.4);
\fill [red] (P)++(0,2.45)circle(0.02);
\fill [red] (P)++(0,2.35)circle(0.02);
\fill [red] (P)++(0,2.25)circle(0.02);

\draw[red] (P) ++(0.4,3) to[out=-120,in=0]++(-0.05,-0.05)to[out=-100,in=30]++(-0.35,-0.85);
\draw[red, densely dotted] (P) ++ (0.4,3) to[out=120, in =60] ++ (-0.3,0);

\draw (P) ++(-0.5,0)++ (-1.5,3) circle[x radius =0.4, y radius=1];
\draw (P) ++(-0.5,0)++ (-1.5,3) ++ (-0.1,0.3) to[out=0, in=0] ++(0,-0.6);
\draw (P) ++(-0.5,0)++ (-1.5,3) ++ (0,0.23) to[out=210,in=150] ++ (0,-0.46);

\draw[red] (P) ++(-0.5,0)++ (-1.12,3) to[out=-120, in =-60] ++ (-0.3,0);
\draw[red, densely dotted] (P) ++(-0.5,0)++ (-1.12,3) to[out=120, in =60] ++ (-0.3,0);

\draw (P) ++(-0.3,0)++ (-1.7,1.5) node{$[\partial D_{i}^{triv}]$};

\draw[->](P)++(-0.3,0)++(-1.1,0) --++(0.7,0);

\draw (P)++(-0.5,0)++(-2.5,1.5) node[left]{$\pi_{i}$};
\draw[->] (P) ++(-0.5,0)++ (-2.5,2)--++(0,-1);

\end{tikzpicture}
\caption{Attaching the trivializing disk by adding annuli}
\label{fig:annulus}
\end{figure}

At first, take the section of $S^{1}$-bundle $\pi_i: E_{i} \rightarrow S_{i}$ restricted to $S_{i}''$ and identify it with $S_{i}''$.
This intersects with the $1$-cycle $t_{i}$ transversally at once and does not intersect with cycles of other types.
However, $S_{i}''$ has a non-empty boundary, thus, it is not a closed cycle.
In general, the $S^{1}$-bundle $\pi_{i}$ may have non-zero Euler number, thus the section $S_{i}''$ cannot extend to $D_{i}^{triv}$.
The boundary of the trivializing disk $D^{triv}_{i}$ is attached to a cycle whose homology class is $h_{i} + (1-|\cP_{i}|)t_{i}$ in $H_{1} (\pi^{-1}_{i}(S_{i}''); \bZ)$, where $h_{i}$ is represented by $\partial D_{i}^{triv} \subset S_{i}''$.
Here, consider $|\cP_{i}|-1$ arcs $\alpha (i; triv, I)$ ($P_{I} \in \cP_{i}$, $I \neq I(i,0)$) and their fiber annuli $A(i; triv, I)$.
We give each annulus $A(i; triv, I)$ an orientation such that $\partial_{I} (A; triv, I)$ is oriented negatively in the plumbing torus $\partial D_{i}^{I} \times \partial D_{I}^{i}$.
Then, the intersection of the union $S_{i}'' \cup \bigcup _{P_{I} \in \cP_{i} \setminus \{P_{I(i,0)}\}} A(i;triv, I)$ with $\pi^{-1}_{i} (D_{i}^{triv})$ represents the homology cycle $h_{i} + (1-|\cP_{i}|)t_{i}$. 
Therefore, the trivializing disk $D_{i}^{triv}$ is attached to the union $S_{i}'' \cup \bigcup_{P_{I} \in \cP_{i} \setminus \{P_{I(i,0)}\} } A(i;triv, I)$ (see Figure \ref{fig:annulus}).
Adding annuli gives a corner to the surface, but suppose they are smoothed (see also the argument and figures in Section 4.1.4 in \cite{doi-hor}).
At this point, $\partial D_{i}^{I(i,0)} \cup \bigcup_{P_{I} \in \cP_{i} \setminus \{P_{I(i,0)}\}} (\partial D_{i}^{I} \cup \partial_{I} A(i;triv, I))$ remains the boundary of the resulting surface.
This operation of attaching the trivializing disk by adding annuli will appear several times after this.

\underline{(ii) The vertex $P_{I(i,0)}$}

Here, take the annulus $A(I(i,0);i,0) \subset \pi^{-1}_{I(i,0)} (S_{I(i,0)}'')$. 
This annulus is attached to $\partial D_{i}^{I(i,0)}$ in the surface constructed in (i) along $\partial_{I(i,0)}A(I(i,0);i,0)$ by plumbing.

\underline{(iii) Vertices $P_{I}$ ($P_{I} \in \cP_{i} \setminus \{P_{I(i,0)}\}$)}

This vertex corresponds to the $S^1$-bundle $\pi_I : E_I \rightarrow S_I$ with the Euler number $-1$.
Similarly as before, take the section restricted to $S''_{I}$ and by extend to the trivializing disk by adding the annulus $A(I; triv, i)$ (which is negatively oriented to the fiber).
Then, the piece $S''_{I} \cup D^{triv}_I \cup A(I;triv, i)$ is attached to $\partial D_{i}^{I} \cup \partial_{I} A(i;triv, I)$ which is a part of the boundary of the surface constructed in (i) by plumbing.
The remaining boundary of the surface is the boundary of plumbing disks $\bigcup_{\ell_{j} \in \ell_{I} \setminus \{\ell_{i}\} } \partial D^{j}_{I}$.

\underline{(iv) Vertices $\ell_{j}$ ($\ell_{j} \notin P_{I(i,0)}$)}

On this vertex, take the annulus $A(j; I(i,j), I(j,0))$ in $\pi^{-1}_{j} (S''_{I})$. Note that $I(i,j) \neq I(j,0)$ since $\ell_{j} \notin P_{I(i,0)}$.
This annulus is attached to $\partial D_{I(i,j)}^{j}$ which is a part of the boundary of the surface constructed in (iii) along $\partial_{I(i,j)} A(j;I(i,j), I(j,0))$ by plumbing.
The other boundary $\partial_{I(j,0)} A(j;I(i,j),I(j,0))$ remains as the boundary of the surface.

\underline{(v) Vertices $P_{I}$ ($P_{I} \in \cP_{0} \setminus \{P_{I(i,0)}\} $)}

This vertex also corresponds to the $S^1$-bundle with the Euler number $-1$ as (iii).
Here, take the section $S''_{I} \cup D_{I}^{triv} \cup A(I;triv, 0)$ extended to the trivializing disk by adding the annulus.
This is attached to $\partial_{I(j,0)} A(j; I(i,j). I(j,0))$ along $\partial D_{I}^{j}$ by plumbing. 
The remaining boundary is $\partial_{0}A(I;triv,0)\cup \partial D_{I}^{0}$.

\underline{(vi) The vertex $\ell_0$}

Here, we construct a surface similar to (i).
By adding $|\cP_{0}| -1 $ annuli $A(0;triv, I)$ ($P_I \in \cP_0 \setminus \{P_{I(i,0)}\}$) whose fiber is negatively oriented to $S_{0}''$ and extend the section to the trivializing disk $D_{0}^{triv}$.
The piece $S''_{0} \cup D_{0}^{triv} \cup \bigcup_{P_{I} \in \cP_{0} \setminus \{P_{I(i,0)}\}} A(0;triv, I)$ has the boundary $\partial D_{0}^{I(i,0)} \cup \bigcup_{P_{I} \in \cP_{0} \setminus \{P_{I}\}} \left ( \partial D_{0}^{I} \cup \partial_{I} A(I;triv,I) \right )$.
The piece $\partial D_{0}^{I(i,0)}$ is attached to $\partial_{0}A(I(i,0);i,0)$ which appeared in (ii) by plumbing.
The other pieces $\partial D_{0}^{I} \cup \partial_{I}A(0;triv,I)$ is attached to $\partial_{0}A(I;triv,0) \cup \partial D_{I}^{0}$ by plumbing.
Finally, we get a surface without boundary. Let $F_{i}$ be the resulting surface.

\begin{proposition}
The surface $F_{i}$ is the intersection dual of $t_{i}$.
\end{proposition}

\proof
At first, observe the intersection with $t$-type cycles.
From (i) of the construction of $F_{i}$, we have $F_{i} \cap t_{i}=+1$.
For $j$ with $j \neq i$ and $\ell_{j} \in P_{I(i,0)}$, the bundles $\pi_{j}:E_{j} \rightarrow S_{j}$ do not appear in the construction of $F_{i}$.
Thus $F_{i}$ do not intersect with $t_{j}$ for such $j$.
When $j$ satisfies $\ell_{j} \notin P_{I(i,0)}$, the intersection product of $F_{i}$ and $t_{j}$ is zero by (iv) (if necessary, move the arc $\alpha(j;I(i,j), I(j,0))$ in order not to pass through the point $x_{j}$).

Next, consider the intersection with $\gamma$-type cycles.
A $\gamma$-type cycle is obtained as the lift of a cycle in the graph $\Gamma(\cC)$, consisting of $\alpha(0; I(j,0),I(j,k))$, $\alpha(I(j,0); 0,j)$, $\alpha(j;I(j,0),I(j,k))$, $\alpha(I(j,k);,j,k)$, $\alpha(k;I(j,k),I(k,0))$ and $\alpha(I(k,0);k,0)$ for some $j$ and $k$.
By the assumption (A2) for the arcs and (vi) in the construction, the piece $\alpha(0;I(j,0),I(k,0))$ do not intersect with $F_{i}$.
By the assumption (A3) and (ii) and (iv) in the construction, the piece $\alpha(I(j,0);0,j)$, $\alpha(I(j,k);j,k)$ and $\alpha(I(k,0);k,0)$ do not intersect with $F_{i}$.
Finally, arcs $\alpha(j;I(j,0),I(j,k))$ and $\alpha(k;I(j,k),I(k,0))$ do not intersect with $F_{i}$ by the assumption (A1) and (iii) in the construction.
(if necessary, perturb the section $S_{i}''$ and $S_{I}''$ all in the above observation). 
Thus, the intersection product of $F_{i}$ with any $\gamma$-type cycle vanishes.

Therefore, the surface $F_{i}$ is the intersection dual of $t_{i}$.
\endproof

\section{Computation of the cohomology ring}
Now, we have an explicit basis for the first and second homology groups of $M$.
Using them, we will compute the cohomology ring.
To compute the cohomology ring, we will compute the intersection product in the homology groups defined via Poincar\'e duality.
Since $H_{k}(M;\bZ)$ is a free module, $H^{k}(M; \bZ) \cong \Hom(H_{k}(M), \bZ)$ and we can take the dual basis of $H_{k}(M;\bZ)$ for each $k$.
For cycles $\alpha, \beta$, we have that $PD(\alpha) \cup PD(\beta) = PD (\alpha \cap \beta)$ (see Theorem 11.9 in \cite{bre}), where $PD(-)$ denotes the Poincar\'e dual of a cycle.
Thus, to compute the cup product, it suffices to compute the intersection product of the dual cycles. 
Now, since the manifold $M$ is an oriented closed $3$-manifold, only the non-trivial cup product is $\cup: H^1 (M; \bZ) \times H^{1}(M;\bZ) \rightarrow H^{2}(M;\bZ)$. Thus we will compute the intersection product $\cap: H_2 (M;\bZ) \times H_{2} (M;\bZ) \rightarrow H_{1}(M;\bZ)$.

\begin{proposition}
The intersection product of $\gamma$-type cycles vanishes. That is, for each $\tau_{j_1, k_1}$ and $\tau_{j_2, k_2}$, $\tau_{j_1, k_1} \cap \tau_{j_2, k_2} = 0$.
\end{proposition}

\proof
The intersection dual of a $\gamma$-type cycle is a plumbing torus.
For distinct edges in a plumbing graph, their corresponding plumbing tori do not intersect. 
Thus, the intersection product of $\gamma$-type cycles vanishes.
\endproof

\begin{proposition}
For the intersection duals $F_{i}$ of $t_{i}$ ($i=1,\ldots, n$) and $\tau_{j,k}$ of $\gamma_{j,k}$ ($(j,k) \in \textbf{nbc} (\cC)$), the intersection product is computed as follows:
\begin{eqnarray*}
F_{i} \cap \tau_{j,k}= 
\left \{
\begin{array}{ll}
-t_{i} + \sum_{m \in I(j,i)} t_{m} & (i=k), \\
-t_{k} & (i \in I(j,k), i \neq k), \\
0 & (\mbox{otherwise}).
\end{array}
\right .
\end{eqnarray*}
\end{proposition}

\begin{figure}[htbp]
\centering
\begin{tikzpicture}
\coordinate (P) at (0,0);
\coordinate (Q) at (2,0);
\coordinate (R) at (2,3);
\coordinate (S) at (5.5,3);

\draw (P) --++(-4,0);
\draw (P)++(0,-2) --++(-4,0);
\draw (P) ++ (0,-1) circle [x radius = 0.5, y radius=1];
\draw (P) ++ (-2,-2) node[below]{$S_{i}''$ (with coorinates $(x,y)$)};

\draw[->] (P) ++ (-3.7, -1) --++ (1,0)node[right]{$x$};
\draw[->] (P) ++ (-3.5,-1.2) --++(0,1.5)node[above]{$z$};
\draw[->] (P) ++ (-3.6,-1.1) --++(0.7,0.7)node[above]{$y$};

\draw[thick, ->, red] (P) ++ (-1,-0.5) --++(0,0.7);
\draw (P)++(-1,0.2)node[above]{normal vector };

\draw (Q) ++ (0,-1) circle [x radius = 0.5, y radius=1];
\draw (Q) ++ (0,-2) node[below]{$\partial D_{i}^{I(j,i)}$};
\draw[->] (Q) ++(0,-1) ++(-0.45,0.4) --++ (0.4,0)node[right]{$x$};
\draw[->] (Q) ++(0,-1) ++(-0.45,0.4) --++ (0.15,0.6)node[above]{$y$};

\draw (R)  circle [x radius = 0.7, y radius=1.2];
\draw (R) ++(0.1,0.5) to[out = 240,in=120] ++(0,-1);
\draw (R) ++(0,0.3) to[out = -60,in=60] ++(0,-0.6);
\draw[->] (R)++(-0.6,0) to[out=-30,in=210] ++(0.5,0)node[below]{$y$};
\draw[->] (R)++(-0.5,-0.15) --++(0,0.5)node[above]{$z$};
\draw (R) ++ (0,-1.2) node[below]{$\pi^{-1}_{i}(\partial D_{i}^{I(j,i)})$};
\draw[->] (R) ++ (0,-2) --++(0,-0.5);

\draw (R) ++ (1.7,0) node{$\supset$};

\draw (S) ++ (0,-1.2) node[below]{$\partial_{I(j,i)} (A;triv, I(j,i)) $};
\draw (S) ++ (0.3,-1.7) node[below]{$\cup \partial D_{i}^{I(j,i)} $ (smoothed)};
\draw (S) ++ (-0.7,-0.1) to[out=90,in=-90] ++(-0.3,0.1)
to[out=90,in=180] ++ (0.3,0.2)
to[out=0,in=90] ++ (0.3,-0.1)
to[out=-90,in=-90] ++ (-0.3,-0); 
\draw[ preaction={line width=3pt, draw=white}] (S) ++ (-0.7,0.1) to[out=90,in=180] ++(0.7,1.1) 
to[out=0,in=90] ++ (0.7,-1.2)
to[out=-90,in=0] ++ (-0.7,-1.2)
to[out=180,in=-90] ++ (-0.7,1.1);
\draw[thick, ->, red] (S) ++ (-0.5,0) --++(0,0.3);
\draw[thick, ->, red] (S) ++ (-0.65,0.6) --++(0.3,0);
\filldraw (S) ++(-0.65,0.6)circle (0.06) node[left]{$p$};

\end{tikzpicture}
\caption{Orientations of normal vectors of $F_{i}$.}
\label{fig:normal}
\end{figure}
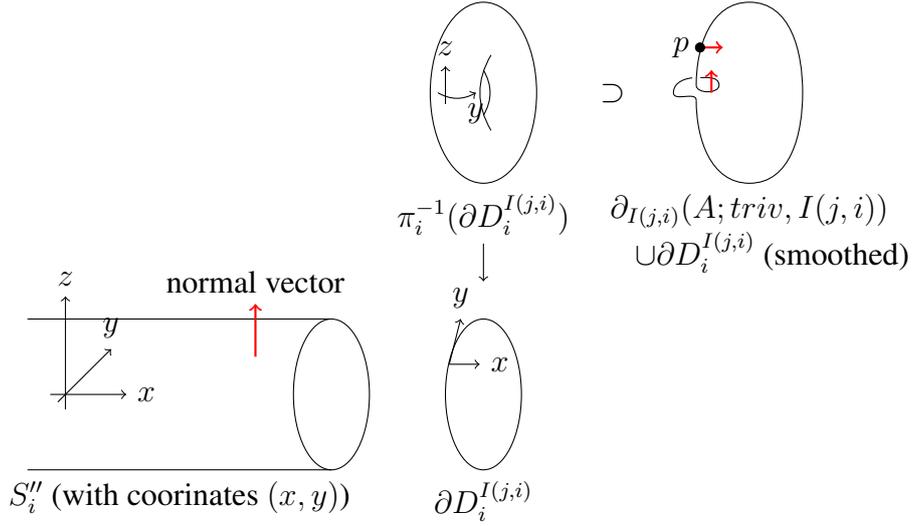

\proof
Suppose that $i=k$. Then, the $F_{i}$ surface and the plumbing torus $\tau_{j,i}$ have a non-empty intersection $\partial_{I(j,i)}A(i;triv, I(j,i)) \cup \partial D_{i}^{I(j,i)}$ by (i) in the construction.
Since the $S^1$-direction of the annulus is negatively oriented, the homology class of the intersection is $\pm(t_{i} - t_{I(j,i)}) = \pm(t_{i} - \sum_{m \in I(j,i)} t_{m} )$ (recall the proof of Proposition \ref{prop:1sthomology}).
To compute the signature, we consider the orientation.
For $p \in \partial_{I(j,i)} A(i; triv, I(j,i)) \subset F_i \cap \tau_{j,i}$, the normal vector of $F_{i}$ points to the positive direction of $\partial D_{i}^{I(j,i)}$ (see Figure \ref{fig:normal}) and of $\tau_{j,k}$ points outward from $S_{i}''$.
Therefore, the intersection is oriented to point in the negative direction of $t_{i}$, thus, the intersection product $F_{i} \cap \tau_{j,i}$ is $- t_{i} + \sum_{m \in I(j,i)} t_{m} $ (recall Convention \ref{conv:orientation}).

Next, suppose that $i \in I(j,k)$ and $i \neq k$.
In this case, the intersection of $F_{i}$ and $\tau_{j,k}$ is $\partial_{I(j,k)}A(k;I(j,k),I(k,0))$ whose homology class is $\pm t_{k}$ by (iii) in the construction (note that $I(j,k)=I(i,j)$).
Again, by a similar observation to the orientation, we have that the intersection product is $- t_{k}$.

Finally, consider the remaining case (that is, $i \notin I(j,k)$).
By the construction of $F_{i}$, a plumbing torus $\tau_{I(j,k)}$ intersects with $F_{i}$ if and only if $i \in I(j,k)$ for $(j,k) \in \textbf{nbc}(\cC)$.
Thus, if $i \notin I(j,k)$, the intersection product $F_{i} \cap \tau_{j,k}$ vanishes.
\endproof

For the intersection duals $F_{i}$ of $t_{i}$ ($i=1,\ldots, n$), we will divide the computation into two cases.

\begin{proposition}\label{prop:intersectionF1}
Suppose that $I(i,0) = I(j,0)$. Then the intersection product $F_{i} \cap F_{j}$ is zero. 
\end{proposition}

\proof
We will observe how $F_{j}$ intersects with $F_{i}$ as $F_{i}$ is constructed through the steps from (i) to (vi). 

Since $\ell_{i} \in P_{I(j,0)}$, the piece corresponding to the vertex $\ell_i$ does not appear in the construction of $F_{j}$. Thus, $F_{i}$ do not intersect with $F_{j}$ around the vertex $\ell_{i}$.

Consider the vertices $P_{I} \in \cP_{i}$.
Around the vertex $P_{I(i,0)}$, the pieces $A(I(i,0);i,0)$ from $F_{i}$ and $A(I(j,0);j,0)$ appear.
By assumption (A4) for arcs, these annuli do not intersect.
If $P_{I} \in \cP_{i}$ satisfies $I \neq I(i,0)$, then $P_{I} \notin \cP_{j}$, since $\cP_{i} \cap \cP_{j} = P_{I(i,0)}=P_{I(j,0)}$. 
Therefore, around such a vertex $P_{I}$, there are no corresponding pieces from $F_{j}$, thus $F_{i}$ does not intersect with $F_{j}$.

Around the vertex $\ell_{k} \, (\notin P_{I(i,0)})$, the pieces $A(k;I(i,k),I(k,0))$ from $F_{i}$ and $A(k;I(j,k),I(k,0))$ appear.
By assumption (A1) for arcs, these pieces do not intersect.

Around the vertex $P_{I} \in \cP_{0}$ ($I \neq I(i,0)$), the piece $S_{I}'' \cup D_{I}^{triv}\cup A(I;triv,0)$ corresponds from both $F_{i}$ and $F_{j}$.
By perturbing the section $S_{I}''$, one allows to separate these pieces.

Finally, the same piece $S_{0}'' \cup D_{0}^{triv} \cup \bigcup_{P_{I} \in P_{0} \setminus \{P_{I(i,0)}\}} A(0;triv,I)$  corresponds to the vertex $\ell_{0}$ from both $F_{i}$ and $F_{j}$.
Similarly, we can separate these pieces around this vertex.
In conclusion, the intersection product $F_{i} \cap F_{j}$ vanishes in this case.

\endproof

\begin{figure}[htbp]
\centering
\begin{tikzpicture}[scale=0.8]
\coordinate (P) at (-0.5,0);
\coordinate (Q1) at (-4,-3);
\coordinate (Q2) at (-6.5,-4);
\coordinate (Q3) at (-4,-4);
\coordinate (Q4) at (-1.5,-4);
\coordinate (Q5) at (-6.5,-7);
\coordinate (Q6) at (-4,-7);
\coordinate (Q7) at (-1.5,-7);
\coordinate (Q8) at (-4,-10);
\coordinate (R1) at (4,-3);
\coordinate (R4) at (6.5,-4);
\coordinate (R2) at (4,-4);
\coordinate (R3) at (1.5,-4);
\coordinate (R7) at (6.5,-7);
\coordinate (R5) at (4,-7);
\coordinate (R6) at (1.5,-7);
\coordinate (R8) at (4,-10);

\coordinate (S1) at (0,-15);
\coordinate (S2) at (-2.5,-16);
\coordinate (S3) at (0,-16);
\coordinate (S4) at (2.5,-16);
\coordinate (S5) at (-2.5,-19);
\coordinate (S6) at (0,-19);
\coordinate (S7) at (2.5,-19);
\coordinate (S8) at (-0,-22);

\draw (P)++(0,-0.3) --++(0,2.3) node[above]{$\ell_1$};
\draw (P)++(-0.2,-0.2) --++(1.7,1.7) node[right]{$\ell_2$};
\draw (P)++(-0.3,0) --++(2.3,0) node[right]{$\ell_3$};
\draw (P)++(-0.3,2)--++(2.3,-2.3)node[below]{$\ell_{0}$} ;


\draw[red] (Q1) ++ (0.5,0) to[out=-90,in=0]++(-0.5,-0.25) to[out=180,in=-90]++(-0.5,0.25);
\draw[red] (Q1) ++ (3,0) to[out=-90,in=0]++(-0.5,-0.25) to[out=180,in=-90]++(-0.5,0.25);
\draw[red] (Q1) ++ (-2,0) to[out=-90,in=0]++(-0.5,-0.25) to[out=180,in=-90]++(-0.5,0.25);
\draw[densely dotted,red] (Q1) ++ (0.5,0) to[out=90,in=0] ++(-0.5,0.25)to[out=180,in=90] ++(-0.5,-0.25);
\draw[densely dotted,red] (Q1) ++ (3,0) to[out=90,in=0] ++(-0.5,0.25)to[out=180,in=90] ++(-0.5,-0.25);
\draw[densely dotted,red] (Q1) ++ (-2,0) to[out=90,in=0] ++(-0.5,0.25)to[out=180,in=90] ++(-0.5,-0.25);
\draw (Q1) ++ (0.5,0)to[out=90,in=180]++(0.75,0.5)to[out=0,in=90]++(0.75,-0.5);
\draw (Q1) ++ (-2,0)to[out=90,in=180]++(0.75,0.5)to[out=0,in=90]++(0.75,-0.5);
\draw (Q1) ++ (-3,0)to[out=90,in=180]++(3,2)to[out=0,in=90]++(3,-2);
\draw (Q1) ++ (0,1.5) circle[radius=0.15];
\draw[blue] (Q1) ++(0,1.35)--++(0,-1.6);
\draw[blue] (Q1) ++(0.15,1.5)to[out=0,in=90]++(2.5,-1.75);
\draw (Q1) ++(2,2)node{$\ell_0$};

\draw (Q2) circle[x radius=0.5, y radius=0.25];
\draw (Q2) ++ (0.5,0)--++(0,-2);
\draw (Q2) ++ (-0.5,0)--++(0,-2);
\draw (Q2) ++ (0.5,-2) to[out=-90,in=0] ++(-0.5,-0.25)to[out=180,in=-90] ++(-0.5,0.25);
\draw[densely dotted] (Q2) ++ (0.5,-2) to[out=90,in=0] ++(-0.5,0.25)to[out=180,in=90] ++(-0.5,-0.25);
\draw(Q2)++(0.5,-1)node[right]{$P_{01}$};
\draw[red] (Q2)++(0,-0.25)--++(0,-2);

\draw[blue] (Q3) circle[x radius=0.5, y radius=0.25];
\draw (Q3) ++ (0.5,0)--++(0,-2);
\draw (Q3) ++ (-0.5,0)--++(0,-2);
\draw[blue] (Q3) ++ (0.5,-2) to[out=-90,in=0] ++(-0.5,-0.25)to[out=180,in=-90] ++(-0.5,0.25);
\draw[densely dotted,blue] (Q3) ++ (0.5,-2) to[out=90,in=0] ++(-0.5,0.25)to[out=180,in=90] ++(-0.5,-0.25);
\draw(Q3)++(0.5,-1)node[right]{$P_{02}$};
\draw (Q3)++(-0.25,-1)circle[radius=0.15];
\draw[red] (Q3)++(0,-0.25)to[out=-90,in=0]++(-0.1,-0.75);

\draw[blue] (Q4) circle[x radius=0.5, y radius=0.25];
\draw (Q4) ++ (0.5,0)--++(0,-2);
\draw (Q4) ++ (-0.5,0)--++(0,-2);
\draw[blue] (Q4) ++ (0.5,-2) to[out=-90,in=0] ++(-0.5,-0.25)to[out=180,in=-90] ++(-0.5,0.25);
\draw[densely dotted,blue] (Q4) ++ (0.5,-2) to[out=90,in=0] ++(-0.5,0.25)to[out=180,in=90] ++(-0.5,-0.25);
\draw(Q4)++(0.5,-1)node[right]{$P_{03}$};
\draw (Q4)++(-0.25,-1)circle[radius=0.15];
\draw[red] (Q4)++(0,-0.25)to[out=-90,in=0]++(-0.1,-0.75);

\draw[red] (Q5) circle[x radius=0.5, y radius=0.25];
\draw (Q5) ++ (0.5,0)--++(0,-2);
\draw (Q5) ++ (-0.5,0)--++(0,-2);
\draw[red] (Q5) ++ (0.5,-2) to[out=-90,in=0] ++(-0.5,-0.25)to[out=180,in=-90] ++(-0.5,0.25);
\draw[densely dotted,red] (Q5) ++ (0.5,-2) to[out=90,in=0] ++(-0.5,0.25)to[out=180,in=90] ++(-0.5,-0.25);
\draw(Q5)++(0.5,-1)node[right]{$\ell_{1}$};
\draw (Q5)++(-0.25,-1)circle[radius=0.15];
\draw[blue] (Q5)++(0,-2.25)to[out=90,in=0]++(-0.1,1.25);

\draw (Q6) circle[x radius=0.5, y radius=0.25];
\draw (Q6) ++ (0.5,0)--++(0,-2);
\draw (Q6) ++ (-0.5,0)--++(0,-2);
\draw (Q6) ++ (0.5,-2) to[out=-90,in=0] ++(-0.5,-0.25)to[out=180,in=-90] ++(-0.5,0.25);
\draw[densely dotted] (Q6) ++ (0.5,-2) to[out=90,in=0] ++(-0.5,0.25)to[out=180,in=90] ++(-0.5,-0.25);
\draw(Q6)++(0.5,-1)node[right]{$\ell_{2}$};
\draw[blue] (Q6)++(0,-0.25)--++(0,-2);

\draw (Q7) circle[x radius=0.5, y radius=0.25];
\draw (Q7) ++ (0.5,0)--++(0,-2);
\draw (Q7) ++ (-0.5,0)--++(0,-2);
\draw (Q7) ++ (0.5,-2) to[out=-90,in=0] ++(-0.5,-0.25)to[out=180,in=-90] ++(-0.5,0.25);
\draw[densely dotted] (Q7) ++ (0.5,-2) to[out=90,in=0] ++(-0.5,0.25)to[out=180,in=90] ++(-0.5,-0.25);
\draw(Q7)++(0.5,-1)node[right]{$\ell_{3}$};
\draw[blue] (Q7)++(0,-0.25)--++(0,-2);

\draw[blue] (Q8) circle[x radius=0.5, y radius=0.25];
\draw[blue] (Q8)++(2.5,0) circle[x radius=0.5, y radius=0.25];
\draw[blue] (Q8)++(-2.5,0) circle[x radius=0.5, y radius=0.25];
\draw (Q8)++(0.5,0)to[out=-90,in=180] ++(0.75,-0.5)to[out=0,in=-90] ++(0.75,0.5);
\draw (Q8)++(-0.5,0)to[out=-90,in=0] ++(-0.75,-0.5)to[out=180,in=-90] ++(-0.75,0.5);
\draw (Q8) ++ (3,0) to[out=-90,in=0] ++(-3,-2);
\draw (Q8) ++ (-3,0) to[out=-90,in=180] ++(3,-2);
\draw (Q8) ++ (0,-1.5) circle[radius=0.15];
\draw[red] (Q8)++(-0.15,-1.5) to[out=180,in=-90]++(-2.35,1.25);
\draw (Q8) ++(2,-2)node{$P_{123}$};
\draw (Q8) ++(0,-2.5)node{$ID(t_{1})$};


\draw[red] (R1) ++ (0.5,0) to[out=-90,in=0]++(-0.5,-0.25) to[out=180,in=-90]++(-0.5,0.25);
\draw[red] (R1) ++ (3,0) to[out=-90,in=0]++(-0.5,-0.25) to[out=180,in=-90]++(-0.5,0.25);
\draw[red] (R1) ++ (-2,0) to[out=-90,in=0]++(-0.5,-0.25) to[out=180,in=-90]++(-0.5,0.25);
\draw[densely dotted,red] (R1) ++ (0.5,0) to[out=90,in=0] ++(-0.5,0.25)to[out=180,in=90] ++(-0.5,-0.25);
\draw[densely dotted,red] (R1) ++ (3,0) to[out=90,in=0] ++(-0.5,0.25)to[out=180,in=90] ++(-0.5,-0.25);
\draw[densely dotted,red] (R1) ++ (-2,0) to[out=90,in=0] ++(-0.5,0.25)to[out=180,in=90] ++(-0.5,-0.25);
\draw (R1) ++ (0.5,0)to[out=90,in=180]++(0.75,0.5)to[out=0,in=90]++(0.75,-0.5);
\draw (R1) ++ (-2,0)to[out=90,in=180]++(0.75,0.5)to[out=0,in=90]++(0.75,-0.5);
\draw (R1) ++ (-3,0)to[out=90,in=180]++(3,2)to[out=0,in=90]++(3,-2);
\draw (R1) ++ (0,1.5) circle[radius=0.15];
\draw[blue] (R1) ++(-0.15,1.5)to[out=180,in=90]++(-2.5,-1.75);
\draw[blue] (R1) ++(0.15,1.5)to[out=0,in=90]++(2.5,-1.75);
\draw (R1) ++(2,2)node{$\ell_0$};

\draw (R2) circle[x radius=0.5, y radius=0.25];
\draw (R2) ++ (0.5,0)--++(0,-2);
\draw (R2) ++ (-0.5,0)--++(0,-2);
\draw (R2) ++ (0.5,-2) to[out=-90,in=0] ++(-0.5,-0.25)to[out=180,in=-90] ++(-0.5,0.25);
\draw[densely dotted] (R2) ++ (0.5,-2) to[out=90,in=0] ++(-0.5,0.25)to[out=180,in=90] ++(-0.5,-0.25);
\draw(R2)++(0.5,-1)node[right]{$P_{02}$};
\draw[red] (R2)++(0,-0.25)--++(0,-2);

\draw[blue] (R3) circle[x radius=0.5, y radius=0.25];
\draw (R3) ++ (0.5,0)--++(0,-2);
\draw (R3) ++ (-0.5,0)--++(0,-2);
\draw[blue] (R3) ++ (0.5,-2) to[out=-90,in=0] ++(-0.5,-0.25)to[out=180,in=-90] ++(-0.5,0.25);
\draw[densely dotted,blue] (R3) ++ (0.5,-2) to[out=90,in=0] ++(-0.5,0.25)to[out=180,in=90] ++(-0.5,-0.25);
\draw(R3)++(0.5,-1)node[right]{$P_{01}$};
\draw (R3)++(-0.25,-1)circle[radius=0.15];
\draw[red] (R3)++(0,-0.25)to[out=-90,in=0]++(-0.1,-0.75);

\draw[blue] (R4) circle[x radius=0.5, y radius=0.25];
\draw (R4) ++ (0.5,0)--++(0,-2);
\draw (R4) ++ (-0.5,0)--++(0,-2);
\draw[blue] (R4) ++ (0.5,-2) to[out=-90,in=0] ++(-0.5,-0.25)to[out=180,in=-90] ++(-0.5,0.25);
\draw[densely dotted,blue] (R4) ++ (0.5,-2) to[out=90,in=0] ++(-0.5,0.25)to[out=180,in=90] ++(-0.5,-0.25);
\draw(R4)++(0.5,-1)node[right]{$P_{03}$};
\draw (R4)++(-0.25,-1)circle[radius=0.15];
\draw[red] (R4)++(0,-0.25)to[out=-90,in=0]++(-0.1,-0.75);

\draw[red] (R5) circle[x radius=0.5, y radius=0.25];
\draw (R5) ++ (0.5,0)--++(0,-2);
\draw (R5) ++ (-0.5,0)--++(0,-2);
\draw[red] (R5) ++ (0.5,-2) to[out=-90,in=0] ++(-0.5,-0.25)to[out=180,in=-90] ++(-0.5,0.25);
\draw[densely dotted,red] (R5) ++ (0.5,-2) to[out=90,in=0] ++(-0.5,0.25)to[out=180,in=90] ++(-0.5,-0.25);
\draw(R5)++(0.5,-1)node[right]{$\ell_{2}$};
\draw (R5)++(-0.25,-1)circle[radius=0.15];
\draw[blue] (R5)++(0,-2.25)to[out=90,in=0]++(-0.1,1.25);

\draw (R6) circle[x radius=0.5, y radius=0.25];
\draw (R6) ++ (0.5,0)--++(0,-2);
\draw (R6) ++ (-0.5,0)--++(0,-2);
\draw (R6) ++ (0.5,-2) to[out=-90,in=0] ++(-0.5,-0.25)to[out=180,in=-90] ++(-0.5,0.25);
\draw[densely dotted] (R6) ++ (0.5,-2) to[out=90,in=0] ++(-0.5,0.25)to[out=180,in=90] ++(-0.5,-0.25);
\draw(R6)++(0.5,-1)node[right]{$\ell_{1}$};
\draw[blue] (R6)++(0,-0.25)--++(0,-2);

\draw (R7) circle[x radius=0.5, y radius=0.25];
\draw (R7) ++ (0.5,0)--++(0,-2);
\draw (R7) ++ (-0.5,0)--++(0,-2);
\draw (R7) ++ (0.5,-2) to[out=-90,in=0] ++(-0.5,-0.25)to[out=180,in=-90] ++(-0.5,0.25);
\draw[densely dotted] (R7) ++ (0.5,-2) to[out=90,in=0] ++(-0.5,0.25)to[out=180,in=90] ++(-0.5,-0.25);
\draw(R7)++(0.5,-1)node[right]{$\ell_{3}$};
\draw[blue] (R7)++(0,-0.25)--++(0,-2);

\draw[blue] (R8) circle[x radius=0.5, y radius=0.25];
\draw[blue] (R8)++(2.5,0) circle[x radius=0.5, y radius=0.25];
\draw[blue] (R8)++(-2.5,0) circle[x radius=0.5, y radius=0.25];
\draw (R8)++(0.5,0)to[out=-90,in=180] ++(0.75,-0.5)to[out=0,in=-90] ++(0.75,0.5);
\draw (R8)++(-0.5,0)to[out=-90,in=0] ++(-0.75,-0.5)to[out=180,in=-90] ++(-0.75,0.5);
\draw (R8) ++ (3,0) to[out=-90,in=0] ++(-3,-2);
\draw (R8) ++ (-3,0) to[out=-90,in=180] ++(3,-2);
\draw (R8) ++ (0,-1.5) circle[radius=0.15];
\draw[red] (R8)++(0,-0.25)--++(0,-1.1);
\draw (R8) ++(2,-2)node{$P_{123}$};
\draw (R8) ++(0,-2.5)node{$ID(t_{2})$};


\draw (S1) ++ (0.5,0) to[out=-90,in=0]++(-0.5,-0.25) to[out=180,in=-90]++(-0.5,0.25);
\draw (S1) ++ (3,0) to[out=-90,in=0]++(-0.5,-0.25) to[out=180,in=-90]++(-0.5,0.25);
\draw (S1) ++ (-2,0) to[out=-90,in=0]++(-0.5,-0.25) to[out=180,in=-90]++(-0.5,0.25);
\draw[densely dotted] (S1) ++ (0.5,0) to[out=90,in=0] ++(-0.5,0.25)to[out=180,in=90] ++(-0.5,-0.25);
\draw[densely dotted] (S1) ++ (3,0) to[out=90,in=0] ++(-0.5,0.25)to[out=180,in=90] ++(-0.5,-0.25);
\draw[densely dotted] (S1) ++ (-2,0) to[out=90,in=0] ++(-0.5,0.25)to[out=180,in=90] ++(-0.5,-0.25);
\draw (S1) ++ (0.5,0)to[out=90,in=180]++(0.75,0.5)to[out=0,in=90]++(0.75,-0.5);
\draw (S1) ++ (-2,0)to[out=90,in=180]++(0.75,0.5)to[out=0,in=90]++(0.75,-0.5);
\draw (S1) ++ (-3,0)to[out=90,in=180]++(3,2)to[out=0,in=90]++(3,-2);
\draw[blue] (S1)++(-2.5,-0.25)to[out=90,in=180]++(1.25,1.3)to[out=0,in=90]++(1.25,-1.3);
\draw (S1) ++(2,2)node{$\ell_0$};

\draw (S2) circle[x radius=0.5, y radius=0.25];
\draw (S2) ++ (0.5,0)--++(0,-2);
\draw (S2) ++ (-0.5,0)--++(0,-2);
\draw (S2) ++ (0.5,-2) to[out=-90,in=0] ++(-0.5,-0.25)to[out=180,in=-90] ++(-0.5,0.25);
\draw[densely dotted] (S2) ++ (0.5,-2) to[out=90,in=0] ++(-0.5,0.25)to[out=180,in=90] ++(-0.5,-0.25);
\draw(S2)++(0.5,-1)node[right]{$P_{01}$};
\draw[red] (S2)++(0,-0.25)--++(0,-2);

\draw (S3) circle[x radius=0.5, y radius=0.25];
\draw (S3) ++ (0.5,0)--++(0,-2);
\draw (S3) ++ (-0.5,0)--++(0,-2);
\draw (S3) ++ (0.5,-2) to[out=-90,in=0] ++(-0.5,-0.25)to[out=180,in=-90] ++(-0.5,0.25);
\draw[densely dotted] (S3) ++ (0.5,-2) to[out=90,in=0] ++(-0.5,0.25)to[out=180,in=90] ++(-0.5,-0.25);
\draw(S3)++(0.5,-1)node[right]{$P_{02}$};
\draw [red] (S3)++(0,-0.25)--++(0,-2);

\draw (S4) circle[x radius=0.5, y radius=0.25];
\draw (S4) ++ (0.5,0)--++(0,-2);
\draw (S4) ++ (-0.5,0)--++(0,-2);
\draw (S4) ++ (0.5,-2) to[out=-90,in=0] ++(-0.5,-0.25)to[out=180,in=-90] ++(-0.5,0.25);
\draw[densely dotted] (S4) ++ (0.5,-2) to[out=90,in=0] ++(-0.5,0.25)to[out=180,in=90] ++(-0.5,-0.25);
\draw(S4)++(0.5,-1)node[right]{$P_{03}$};

\draw (S5) circle[x radius=0.5, y radius=0.25];
\draw (S5) ++ (0.5,0)--++(0,-2);
\draw (S5) ++ (-0.5,0)--++(0,-2);
\draw (S5) ++ (0.5,-2) to[out=-90,in=0] ++(-0.5,-0.25)to[out=180,in=-90] ++(-0.5,0.25);
\draw[densely dotted] (S5) ++ (0.5,-2) to[out=90,in=0] ++(-0.5,0.25)to[out=180,in=90] ++(-0.5,-0.25);
\draw(S5)++(0.5,-1)node[right]{$\ell_{1}$};
\draw[blue] (S5)++(0,-0.25)--++(0,-2);

\draw (S6) circle[x radius=0.5, y radius=0.25];
\draw (S6) ++ (0.5,0)--++(0,-2);
\draw (S6) ++ (-0.5,0)--++(0,-2);
\draw (S6) ++ (0.5,-2) to[out=-90,in=0] ++(-0.5,-0.25)to[out=180,in=-90] ++(-0.5,0.25);
\draw[densely dotted] (S6) ++ (0.5,-2) to[out=90,in=0] ++(-0.5,0.25)to[out=180,in=90] ++(-0.5,-0.25);
\draw(S6)++(0.5,-1)node[right]{$\ell_{2}$};
\draw[blue] (S6)++(0,-0.25)--++(0,-2);

\draw (S7) circle[x radius=0.5, y radius=0.25];
\draw (S7) ++ (0.5,0)--++(0,-2);
\draw (S7) ++ (-0.5,0)--++(0,-2);
\draw (S7) ++ (0.5,-2) to[out=-90,in=0] ++(-0.5,-0.25)to[out=180,in=-90] ++(-0.5,0.25);
\draw[densely dotted] (S7) ++ (0.5,-2) to[out=90,in=0] ++(-0.5,0.25)to[out=180,in=90] ++(-0.5,-0.25);
\draw(S7)++(0.5,-1)node[right]{$\ell_{3}$};

\draw (S8) circle[x radius=0.5, y radius=0.25];
\draw (S8)++(2.5,0) circle[x radius=0.5, y radius=0.25];
\draw (S8)++(-2.5,0) circle[x radius=0.5, y radius=0.25];
\draw (S8)++(0.5,0)to[out=-90,in=180] ++(0.75,-0.5)to[out=0,in=-90] ++(0.75,0.5);
\draw (S8)++(-0.5,0)to[out=-90,in=0] ++(-0.75,-0.5)to[out=180,in=-90] ++(-0.75,0.5);
\draw (S8) ++ (3,0) to[out=-90,in=0] ++(-3,-2);
\draw (S8) ++ (-3,0) to[out=-90,in=180] ++(3,-2);
\draw[red] (S8)++(-2.5,-0.25)to[out=-90,in=180]++(1.25,-1)to[out=0,in=-90]++(1.25,1);
\draw (S8) ++(2,-2)node{$P_{123}$};
\draw (S8) ++(0,-2.5)node{$ID(t_{1}) \cdot ID(t_{2}) = \gamma_{1,2}$};


\end{tikzpicture}
\caption{An example of the intersection product}
\label{fig:intersection}
\end{figure}
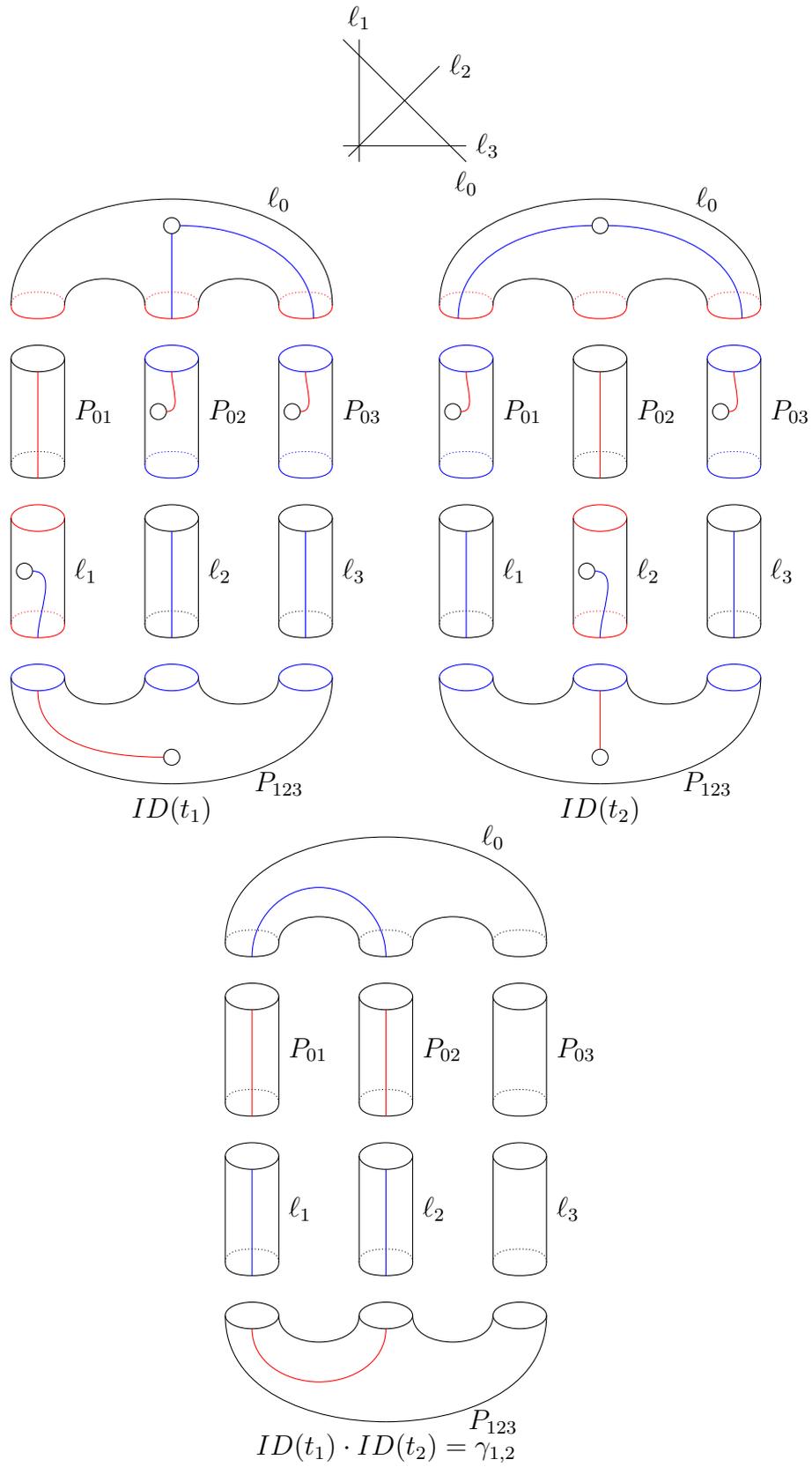

\begin{proposition}\label{prop:intersectionF2}
Suppose that $I(i,0) \neq I(j,0)$. Then the intersection product is computed as follows (see Figure \ref{fig:intersection}):
\begin{eqnarray*}
F_{i} \cap F_{j}= 
\left \{
\begin{array}{ll}
\gamma_{i,j} & ((i,j) \in nbc(\cC) ), \\
\gamma_{k,j} - \gamma_{k,i} & (\exists k \mbox{ s.t. $(k,i),(k,j) \in nbc(\cC)$}).
\end{array}
\right .
\end{eqnarray*}
\end{proposition}

\proof
Similarly to before, we will observe the intersection according to the steps (i) through (vi). 

In this case, recall that there exists a unique index $I(i,j) \, (\neq I(i,0), I(j,0))$ such that $\ell_i , \ell_{j} \in P_{I(i,j)}$.
Since $\ell_{i} \notin P_{I(j,0)}$, 
the pieces $S_{i}'' \cup D_{i}^{triv} \cup \bigcup_{P_{I} \in \cP_{i} \setminus \{P_{I(i,0)}\}}A(i;triv, I)$ from $F_{i}$ (by (i)) and $A(i; I(i,j), I(i,0))$ from $F_{j}$ correspond to the vertex $\ell_{i}$ (by (iv)).
Their intersection is an arc $A(i; I(i,j), I(i,0)) \cap S_{i}''$.
And this, by perturbation, can be expressed as $\alpha (i; I(i,j), I(i,0))$ as defined before.

To the vertex $P_{I(i,0)}$, the piece $A(I(i,j);i,0)$ from $F_{i}$ (by (ii)) and $S_{I(i,0)}''\cup D_{I(i,0)}^{triv} \cup A(I(i,0); triv, 0)$ from $F_{j}$ (by (v)) correspond.
Their intersection is an arc $A(I(i,0); i,0) \cap S_{I(i,0)}''$ and it can be expressed as $\alpha (I(i,0);i,0)$.

Among the vertices $P_{I}$ ($P_{I} \in \cP_{i} \setminus \{P_{I(i,0)}\}$), only for the vertex $P_{I(i,j)}$, the pieces obtained from both $F_{i}$ and $F_{j}$ correspond.
The pieces $S_{I(i,j)}'' \cup D_{I(i,j)}^{triv} \cup A(I(i,j); triv, i)$ from $F_{i}$ and $S_{I(i,j)}'' \cup D_{I(i,j)}^{triv} \cup A(I(i,j); triv, j)$ from $F_{j}$ correspond to the vertex $P_{I(i,j)}$ (both by (ii)).
Their intersection is the union of arcs $\alpha (I(i,j); triv, i) \cup \alpha (I(i,j); triv, j)$, and they are connected in the trivializing disk. 
The resulting arc is expressed as $\alpha (I(i,j);i,j)$ (if necessary, by perturbation), and this arc is the intersection on the vertex $P_{I(i,j)}$.

Consider the vertex $\ell_{k}$ ($\ell_{k} \notin P_{I(i,0)}$).
When $k \neq j$, the piece $A(k;I(j,k), I(k,0))$ is obtained from $F_{j}$ (by (iv)).
By the assumption (A1) for the arcs, this annulus does not intersect with $A(k; I(i,k), I(k,0))$, which is the piece obtained from $F_{i}$ (by (iv)).
When $k=j$, the piece $S_{j}'' \cup D_{j}^{triv} \cup \bigcup_{P_{I} \in \cP_{j} \setminus \{P_{I(j,0)}\}}$ is obtained from $F_{j}$ (by (i)).
This piece intersects with the annulus $A(j; I(i,j), I(j,0))$ from $F_{i}$ along the arc, which can be expressed as $\alpha (j; I(i,j), I(j,0))$.

For the vertices $P_{I} \in \cP_{0}$, when $I \neq I(j,0)$, the piece $S_{I}'' \cup A (I; triv, 0)$ is obtained from both $F_{i}$ and $F_{j}$ (by (v)).
By a small perturbation, we can separate these pieces, and thus they do not intersect.
For the vertex $P_{I(j,0)}$, the intersection of pieces from $F_{i}$ and $F_{j}$ is the arc expressed as $\alpha (I(j,0); j,0)$, from the same reason on the vertex $P_{I(i,0)}$.

Finally, consider the vertex $\ell_{0}$.
Around this vertex, the piece $S_{0}'' \cup D_{0}^{triv} \cup \bigcup_{P_{I} \in \cP_{0} \setminus \{P_{I(i,0)}\}} A(0;triv, I)$ is obtained from $F_{i}$ and $S_{0}'' \cup D_{0}^{triv} \cup \bigcup_{P_{I} \in \cP_{0} \setminus \{P_{I(j,0)}\}} A(0;triv, I)$ is obtained from $F_{j}$ (both from (vi)).
Near the boundary $\partial D_{0}^{I}$ ($I \neq I(j,0), I(k,0)$), both pieces represents the same homology class, namely, $t_{0} - t_{I}$.
Therefore, we can separate the pieces around these boundaries.
However, around the boundary $\partial D_{0}^{I(i,0)}$ and $\partial D_{0}^{I(j,0)}$,
the annulus and the section obtained from each component correspond.
Thus, by perturbation, the intersection is the union $(A (0;triv, I(i,0)) \cap S_{0}'' ) \cup (A(0;triv, I(j,0)) \cap S_{0}'')$.
They are connected in the trivializing disk similarly as discussed before.
Therefore, the resulting intersection can be expressed as $\alpha(0;I(i,0),I(j,0))$.

As a result, the intersection of $F_{i}$ and $F_{j}$ is a cycle concatenating the arcs $\alpha(0;I(i,0),I(j,0))$, $\alpha(I(j,0); 0,j)$, $\alpha(j; I(j,0), I(i,j))$, $\alpha(I(i,j);i,j)$, $\alpha(i; I(i,j),I(i,0))$ and $\alpha(I(i,0);0,i)$.
Finally, we check the orientation of this cycle.
Take $p \in \partial D_{i}^{I(i,0)} \cap \partial_{I(i,0)}A(0;triv, I(i,0))$ for example (actually this set consists of a single point).
The normal vector of $F_{i}$ points to the positive direction of $t_{i}$ and of $F_{j}$ points to the positive direction of $\partial D_{i}^{I(i,0)}$ (note that the fiber direction of the annulus $A(0;triv. I(i,0))$ is negatively oriented) at $p$.
Thus, the tangent vector points outward from $S_{i}''$ at $p$.
This agrees with the orientation of the $\gamma$-type cycle, which is defined just before Remark \ref{rmk:realizable1cycle}.
By the definition of the cycle $\gamma_{i,j}$, the proposition holds.
\endproof

By summarizing Propositions \ref{prop:intersectionF1} and \ref{prop:intersectionF2}, we have the following.

\begin{proposition}
For the intersection duals $F_{i}$ of $t_{i}$ ($i=1,\ldots, n$), the intersection product is computed as follows:
\begin{eqnarray*}
F_{i} \cap F_{j}= 
\left \{
\begin{array}{ll}
\gamma_{i,j} & ((i,j) \in nbc(\cC) ), \\
\gamma_{k,j} - \gamma_{k,i} & (\exists k \mbox{ s.t. $(k,i),(k,j) \in nbc(\cC)$}), \\
0 & (\mbox{otherwise}).
\end{array}
\right .
\end{eqnarray*}
\end{proposition}

The above computation results allow us to compute the cohomology ring of $M$.

\begin{theorem}
The cohomology ring $H^{*} (M(\cC); \bZ)$ is isomorphic to the double of the Orlik-Solomon algebra.
\end{theorem}

\proof
As a basis of the first cohomology group $H^{1} (M(\cC);\bZ)$, take the dual basis $\{\overline{t}_{i} \}$ and $\{\overline{\gamma_{j,k}}\}$ of $\{t_{i} \}$ and $\{ \gamma_{j,k} \}$.
The basis $\{ F_{i} \}$ and $\{ \tau_{j,k} \}$ of the second homology group, which we constructed, are the intersection duals.
Therefore, it follows that $PD(F_{i}) = \overline{t}_{i}$ and $PD(\tau_{j,k}) = \overline{\gamma}_{j,k}$.
Similarly for the second cohomology group, it follows that $PD(t_{i}) = \overline{F}_{i}$ and $PD(\gamma_{j,k}) = \overline{\tau}_{j,k}$.
By defining a homomorphism $H^{*}(M;\bZ) \rightarrow \widehat{A}$ by
$\overline{t}_{i} \mapsto e_{i}, \overline{\gamma}_{j,k} \mapsto \overline{f}_{j,k}, \overline{F}_{i} \mapsto \overline{e}_{i}, \overline{\tau}_{j,k} \mapsto f_{j,k}$,
the theorem follows from the above computation result and $PD(a) \cup PD(b) = PD(a\cap b)$.
\endproof

\section{Resonance varieties}

As an application of the main result, we compute the resonance variety of the boundary manifold of a combinatorial line arrangement.
Let $A=\oplus_{k=0}^{m}A^{k}$ be a connected, graded-commutative graded algebra over a field $\mathbb{K}$. We assume $\mathbb{K}$ is algebraically closed and $\Char \mathbb{K}=0$. We also assume that each $A^{k}$ is finite-dimensional.
Since $a \cdot a =0$ for $a \in A^1$, we have a cochain complex 
\[
(A, a): 0 \xrightarrow{\cdot a} A^{0} \xrightarrow{\cdot a} A^{1} \xrightarrow{\cdot a} \cdots \xrightarrow{\cdot a} A^{\ell} \xrightarrow{\cdot a} 0.
\]
The \textit{resonance variety} of $A$ is defined as its cohomology jumping loci:
\[
\cR_{d}^{k}(A) = \{a \in A^{1} \mid \dim_{\mathbb{K}}H^{k}(A, a) \geq d\}.
\]
The resonance variety is an important object in the study of the topology of hyperplane arrangements, and it has been extensively investigated \cite{fal-yuz, lib-yuz, yuz}.
In \cite{coh-suc-proj}, the resonance variety of doubled algebras was also studied. 
Since the Orlik-Solomon algebra of a combinatorial line arrangement, which we mainly consider, has degree two, we begin by reviewing the resonance variety of the doubling of a graded algebra of degree two.

Let $\{a_{i}\}_{i=1}^{r_{1}}$, $\{b_{j}\}_{j=1}^{r_{2}}$ be basis of $A^1$ and $A^2$, respectively ($r_{k} = \dim_{\bK}A^{k}$).
Suppose that the multiplication $\mu: A^{1} \otimes A^{1} \rightarrow A^{2}$ is described as 
\[
\mu(a_{i}, a_{j}) = \sum_{k=1}^{r_{2}} \mu_{i,j,k}b_{k}.
\]
Let $(a, b) \in A^{1} \oplus \overline{A}^{2}=\mathsf{D}(A)^1$, and assume $a=\sum_{i=1}^{r_1} p_{i} a_{i}$ and $b=\sum_{j=1}^{r_{2}} q_{j} b_{j}$.
Define an $r_{1} \times r_{2}$ matrix $\Delta^{a}$ by 
\[
\Delta^{a}_{j,k} = \sum_{i=1}^{r_{1}} \mu_{i,j,k} p_{i}.
\]
We also define a $r_{1} \times r_{1}$ matrix $\Phi^{b}$ by 
\[
\Phi^{b}_{i,j} = \sum_{k=1}^{r_{2}} \mu_{i,j,k} q_{k}.
\]

The cochain complex defined on the doubling algebra by $(a,b) \in \widehat{A}^1$ is expressed as follows,
\[
(\mathsf{D}(A), (a,b)): A^{0} \xrightarrow{d_{1}} A^{1}\oplus \overline{A}^{2} \xrightarrow{d_{2}} A^{2}\oplus \overline{A}^{1} \xrightarrow{d_{3}} \overline{A}^{0},
\]
where each boundary map is described as
\[ 
d_{1} = 
\begin{pmatrix}
a & b
\end{pmatrix},
\hspace{3mm}
d_{2}= 
\begin{pmatrix}
\Phi^{b} & \Delta^{a} \\
-{}^{t} \! \Delta^{a} & 0
\end{pmatrix},
\hspace{3mm}
d_{3} = 
\begin{pmatrix}
{}^{t}  a \\ {}^{t}  b
\end{pmatrix}.
\]
The resonance variety of the double $\mathsf{D}(A)$ is similarly defined as 
\[
\cR_{d}^{k} (\mathsf{D}(A)) = \{(a,b) \in A^{1} \oplus \overline{A}^2 \mid
\dim_{\bK} H^{k} (\mathsf{D}(A), (a,b)) \geq d\}.
\]

Now, suppose that $A=A\otimes \bC$ is the Orlik-Solomon algebra of a combinatorial line arrangement with complex coefficients. At first, recall that $a \in A^{1}$ is nonresonant if $\dim_{\bK} H^{k}(A,a)$ is minimal, particularly, $H^{0}(A,a) = H^{1}(A,a)=0$ for Orlik-Solomon algebras (see Theorem 7.5 in \cite{yuz}).
We also have $\dim_{\bK}H^{2} (A,a) = 1 - b_{1} (A)+ b_{2} (A)$ for nonresonant $a \in A^1$.
The resonance variety of the boundary manifold is nothing but the resonance variety of the double of the Orlik-Solomon algebra.

\begin{proposition}
Let $\cC$ be a combinatorial line arrangement, $A$ be its Orlik-Solomon algebra, and $\mathsf{D}(A)$ be its double. Let $\beta = 1 - b_{1} (A)+ b_{2} (A)$.
\begin{enumerate}[(1)]
\item For nonresonant $a \in A^1$, $(a,b) \in \mathsf{D}(A)^1$ is nonresonant for any $b \in \overline{A}^2$. For such $(a,b)$, we have $H^{0} (\mathsf{D}(A), (a,b))= H^{3} (\mathsf{D}(A), (a,b)) = 0$ and $H^{1} (\mathsf{D}(A), (a,b)) \cong H^{2} (\mathsf{D}(A), (a,b)) \cong \bK^{\beta}$.
\item If $a \in \cR_{d}^{1} (\mathsf{D}(A))$ and $a\neq 0$, then $(a,b) \in \cR^{1}_{d + \beta} (\mathsf{D}(A))$ for any $b \in \overline{A}^2$.
\item Let $d = b_{2} (A) -1 + \dim_{\bK}(\Ker(\Phi^{b}))$. Then, $(0,b) \in \cR^{1}_{d} (\mathsf{D}(A))$ for any $b \in \overline{A}^2$ with $b \neq 0$.
\end{enumerate}
\end{proposition}

\proof
See Lemma 6.5, Proposition 6.7 in \cite{coh-suc-proj}, and their proof.
Since their proof is purely algebraic, all the arguments therein apply to our setting as well.
\endproof

From this, we deduce the following.

\begin{corollary}\label{cor:resonance}
For the resonance variety of the double $\mathsf{D}(A)$ of the Orlik-Solomon algebra, we have:
\begin{enumerate}[(1)]
\item $\cR_{d}^{1}(\mathsf{D}(A)) = \mathsf{D}(A)^{1}$ for $d \leq \beta$.
\item $\cR_{d}^{1}(A) \times \overline{A}^2 \subset \cR_{d + \beta}^{1} (\mathsf{D}(A))$.
\item $\{0\} \times \cR_{d}(\Phi) \subset \cR^{1}_{d+b_2 (A)} (\mathsf{D}(A))$.
\end{enumerate} 
Here, $\cR_{d}(\Phi)$ is the zero set of the ideal defined by $(b_{1}-d) \times (b_{1}-d)$ minors of the matrix $\Phi^{b}$ regarding $q_{k}$ as variables.
\end{corollary}

Note that the third statement above implies that the resonance variety may contain a non-linear irreducible component (see Proposition 6.10 in \cite{coh-suc-proj}).

Let us consider the special case, $d=1$.
Recall that a combinatorial line arrangement $\cC$ is \textit{pencil} if all the lines are contained in the same intersection point, and \textit{near pencil} if all but one of the lines intersect at a single point.
Such combinatorial line arrangements are realized over $\bC$ by defining equation $(x^{n}-y^{n})$ and $(x^{n-1} - y^{n-1})z$ in $\bC P^2$, for example.
If $\cC$ is either pencil or near pencil, then $\beta = 1 -b_{1} (A) + b_{2} (A) \leq 0$. Note that the converse is true. That is, if $\beta \geq1$, then $\cC$ is neither pencil nor near pencil.
This follows, for example, from an inductive argument on the characteristic polynomial.

Now, the computational results in the case $d=1$ is as follows.

\begin{proposition}
Let $\cC$ be a combinatorial line arrangement with $\cL =\{\ell_0,\ldots, \ell_n\}$ and $\mathsf{D}(A)$ be the double of the Orlik-Solomon algebra. Then, 
\begin{eqnarray*}
\cR_{1}^{1} (\mathsf{D}(A)) \cong
\left \{
\begin{array} {ll}
\bC^{n} & \mbox{($\cC$ is pencil)}, \\
\bC^{2n-2} & \mbox{($\cC$ is near pencil)}, \\
\bC^{b_1(A)+b_{2}(A)} & \mbox{(otherwise)}.
\end{array}
\right .
\end{eqnarray*}
\end{proposition}

\proof
For the case $\cC$ is either pencil or near pencil, see Corollary 8.7 in \cite{coh-suc-bound}. For the other cases, we have $\beta \geq 1$ and thus the result follows from Corollary \ref{cor:resonance} (1).
\endproof

\end{document}